\documentclass[11pt]{article}
\usepackage[pdftex]{graphicx}
\usepackage{vmargin}
\usepackage{amssymb}
\usepackage{fancyheadings}
\usepackage{natbib}
\usepackage{color}
\usepackage[tight]{subfigure}
\usepackage{booktabs}
\usepackage{float}
\usepackage{multirow}
\usepackage[pagewise]{lineno}
\usepackage{setspace}
\usepackage{algorithmic, algorithm}
\usepackage{longtable}
\usepackage{booktabs}
\usepackage{amsmath}
\usepackage{url}

\DeclareMathOperator*{\argmax}{arg\,max}

\setpapersize{USletter}
\setmarginsrb{1.5in}{1in}{1in}{0.5in}{12pt}{11mm}{12pt}{30pt}

\begin{document}
\title{Closed-loop field development with multipoint geostatistics and statistical performance assessment}
\author{Mehrdad G Shirangi\footnote{Stanford University, Stanford, CA 94305, United States. Email addresses: mehr@stanford.edu, Mehrdad.GharibShirangi@ge.com }}

\maketitle

\date{}




\section*{Abstract}
\label{Sec:Abstract} 
Closed-loop field development (CLFD) optimization is a comprehensive framework for optimal development of subsurface resources. CLFD involves three major steps: 1) optimization of full development plan based on current set of models, 2) drilling new wells and collecting new spatial and temporal (production) data, 3) model calibration based on all data.
This process is repeated until the optimal number of wells is drilled.
This work introduces a new CLFD implementation for complex systems described by multipoint geostatistics (MPS). Model calibration is accomplished in two steps: conditioning to spatial data by a geostatistical simulation method, and conditioning to production data by optimization-based PCA.
A statistical procedure (TruMAP) is presented to assess the performance of CLFD.
For performance assessment by TruMAP,
the methodology is applied to an oil reservoir example for 25 different true-model cases.
Application of a single-step of CLFD, improved the true NPV in 64\%--80\% of cases. The full CLFD procedure (with three steps) improved the true NPV in 96\% of cases, with an average improvement of 37\%.
These results indicate the effectiveness of performing multiple steps of closed-loop optimization.
This massive computational experiment involved about 9.5 million reservoir simulation runs that took about 320,000 CPU hours.

\section{Introduction}
For optimal operation of oil and gas resources, engineers
typically build physics-based flow simulation models. These models need to be calibrated
with observed spatial and temporal data so that they can be reliable for purposes such as optimization
and decision making.
The common approach
for management of these systems is closed-loop optimization which
involves three major steps, repeated throughout the project life: 1) optimization of decision parameters,
2) implementation and operation for a time period, 3) and model calibration for consistency with (new) observed data.
For different types of reservoir problems, various approaches have been proposed for Steps 1 and 3, i.e., the optimization and the model calibration step.
In addition, implementations of closed-loop optimization may differ on type of decision parameters included in these two steps.

In the context of conventional oil reservoir management,
decision parameters depend on the recovery process employed.
The most common recovery process is water-flooding where water is injected through injection wells to push the oil towards producer wells.
There are other processes such as gas injection, water-alternating-gas injection, surfactant and polymer flooding, and
steam-assisted gravity drainage (SAGD) which are employed based on subsurface reservoir condition and economic parameters (e.g., availability/cost of surface facilities).
This work focuses on the water-flooding process, though the framework and results can be extended to other recovery processes.
Decision parameters, in this context, include the number, type (oil production/water injection), location, and controls (time-varying well bottom-hole pressures) of new wells
throughout the project life.

Closed-loop field development (CLFD) optimization, introduced by \citet{shirangi:15b},
is a comprehensive reservoir management framework.
In the optimization step of CLFD, the number, type, location and controls
of new wells (together with controls of the existing wells) are optimized.
The model calibration step of CLFD involves integration/assimilation of new
spatial data (hard data from new wells) and temporal data (production data from existing wells).
The impact of CLFD research is significant as drilling new wells is
one of the most expensive parts of reservoir operations.
In \citet{shirangi:15b}, the CLFD framework is presented with two-point geostatistical models.
In recent years, however, there has been significant research on multipoint geostatistics (MPS)
which is able to model complex geological features such as channels and deltaic fans.

In this work, a closed-loop field development framework is developed and applied for
reservoir models described by multipoint geostatistics.
A two-step model calibration procedure is presented where we employ an efficient geostatistical algorithm
to integrate spatial data, and a recent model calibration algorithm for integration
of temporal data (production history).
A statistical assessment procedure is then presented for performance evaluation of this closed-loop framework.
In this assessment, the project outcome (true NPV) is treated as a random variable
and the improvement in true NPV is investigated after each CLFD step.

Before the introduction of CLFD framework, most of the reservoir optimization research
has focused on developing improved approaches for optimal operation of existing wells,
referred to as closed-loop reservoir management (CLRM).
Compared to CLFD, CLRM does not involve drilling new wells.
In optimization step of CLRM,
referred to as production optimization or well control optimization, continuous operational settings of existing
wells are optimized. The model calibration of CLRM only includes new temporal data, and integration of new spatial data
is not typically considered.
CLRM has been investigated extensively (see, e.g., \citet{jansen:09,jansen:05}).
Most papers on CLRM investigated the application of particular model calibration and optimization approaches for water-flooding operations (e.g., \citet{aitokhuehi:05,chen:09ens}).
CLRM has also been applied to SAGD operations \citep{regtien:10},
and for the management of geological carbon storage operations \citep{cameron:14}.
It's worth mentioning that \citet{jansen:05,jansen:09} envisioned the use of optimal position of sidetracks or infill wells in closed-loop optimization,
which is now performed in CLFD.

In the field development planning (FDP) problem, encountered in CLFD, decision parameters may include the number of new wells,
well type (producer or injector), well locations, drilling sequence, and well settings.
Most papers on FDP considered the optimization of only a subset of these parameters, mainly the location
of new wells (referred to as the well placement problem).
Many of the methods developed for this problem entail stochastic (global) search methods such as GAs \citep{ozdogan:06},
PSO \citep{onwunalu:10,arnold:16},
evolution strategy with covariance matrix adaptation (CMA-ES) \citep{bayer:07,bayer:08,bouzarkouna:12},
and differential evolution (DE) \citep{bayer:10,nwankwor:13}.
Local optimization methods such as SPSA \citep{shirangi:12},
and pattern search techniques \citep{wilson:13,cameron:12} have also been applied for well placement optimization.

The more complex problem of jointly optimizing various parameters in FDP has also been investigated.
\citet{bellout:12,li-jafar:12,zhang:17} presented different algorithms for the joint optimization of well location and controls.
\citet{isebor:14a,isebor:14b} developed a formulation based on a hybrid of PSO, a global stochastic search algorithm, and mesh adaptive direct search (MADS), a local pattern search method. This PSO-MADS procedure can simultaneously optimize the number and type (e.g., injector or producer) of new wells and the drilling sequence, in addition to well locations and controls. This algorithm will be used in this work for the optimization step.

Model calibration is the third major step in CLFD.
In model calibration (also referred to as history matching, or data assimilation),
the goal is to learn from data through generating one or more computational models that are consistent with prior
geological information and new spatial data, and provide flow simulation results that closely match observed (temporal) production data.
The model calibration methods typically involve a challenging minimization problem as this problem is
usually ill-posed and the number of unknown model parameters can be very large.
The ill-posedness of model calibration can be mitigated by reducing the number of parameters
through an appropriate parameterization such as TSVD \citep{shirangi:11,shirangi:14,shirangi:16,bjarkason:17,dickstein:17},
ensemble-based methods \citep{rafiee:17,rafiee:17b,rafiee:18}, and PCA \citep{vo:16}.
\citet{vo:15,vo:14} presented a differentiable PCA-based parameterization (O-PCA)
that enables application of efficient gradient-based approaches for model calibration of complex channelized systems.
This optimization-based principal component analysis (O-PCA) approach is used in the CLFD model calibration step in this work.

The CLFD framework has also been applied in recent work by other researchers.
\citet{morosov:16} applied the CLFD to a realistic example,
while \citet{hanea:17} applied the CLFD framework in simpler settings where only drilling sequence is optimized.
Compared with the original work \citep{shirangi:15b,shirangi:13}, they used different reservoir modeling and
different model calibration and optimization approaches in their CLFD application.
Both papers, reported that CLFD significantly improved true-model NPV in most cases;
however, in at least one case in each work, the true-model NPV decreased slightly after the application of CLFD.
More recently, \citet{hidalgo:17} applied the CLFD algorithm to a realistic reservoir case and reported 40.8\% improvement in (true) NPV.


In \citet{shirangi:15b}, reservoir properties were represented in terms of two-point spatial statistics.
The first contribution of this work is to present and apply a CLFD implementation for complex channelized models, described by a training image.
A training image contains a range of valuable information such as orientation, proportions, affinity, trend that need to be reproduced \citep{tahmasebi:14,mariethoz:14}.
This extension involves particular treatments for integration of spatial and temporal data in the model calibration step.
This new implementation is to improve CLFD model calibration step by utilizing a multipoint geostistical algorithm within the closed-loop optimization workflow.
The second contribution of this work is to present and apply the TruMAP (true-model based assessment of performance)
procedure.
The goal here is to statistically characterize the performance of a closed-loop optimization algorithm in terms of improvement in true NPV.


This paper proceeds as follow.
In the next section, the CLFD framework and treatments for the model calibration and the optimization step are described.
In the section after that, computational results for CLFD application are presented for a two-dimensional binary channelized reservoir case.
Results of statistical assessment of CLFD performance is then presented in this section.
Finally, conclusions and suggestions for future work are provided.

\section{Methodology: closed-loop field development optimization} \label{clfd}
\begin{figure}
\centering
    \includegraphics[width=0.8\textwidth]
    {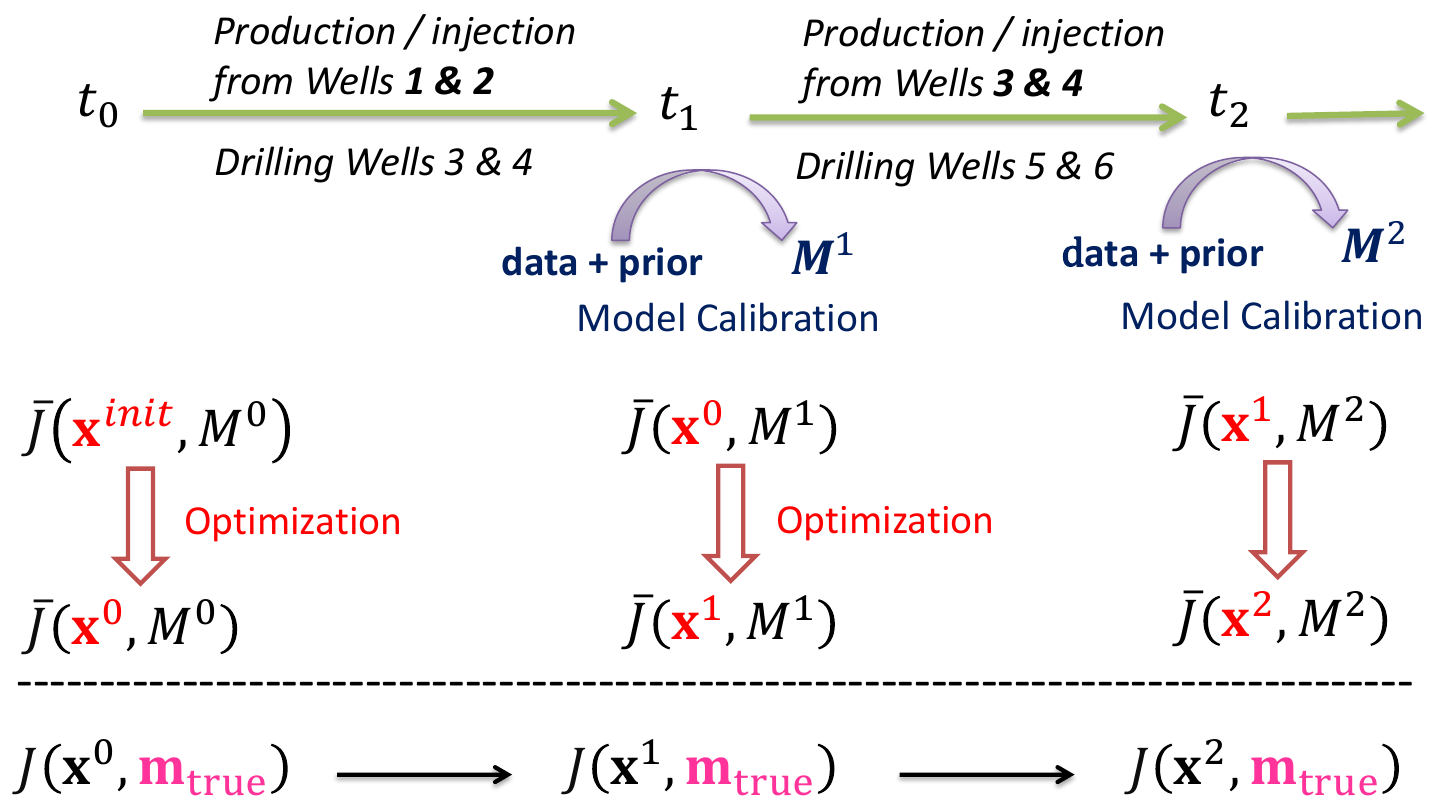}
\caption{Schematic and notations for CLFD framework with $n_{\text{rig}}$=2. }
\label{fig:schem-clfd}
\end{figure}

Closed-loop field development (CLFD) (Fig.~\ref{fig:schem-clfd})
has three major steps as follows.
1) Optimization: optimize the ``full development plan'' (the number, type, location and controls of all new wells together with controls of existing wells),
for a set of models describing the current state of knowledge,
2) operation: drill and complete $n_{\text{rig}}$ new wells, produce from existing wells, and collect new spatial and temporal data,
3) model calibration (history matching): update models for consistency with new data.
These steps are repeated until drilling additional wells does not (sufficiently) increase the expected NPV, at which step, the CLRM workflow is applied.
We now describe the CLFD workflow in more detail.

We consider the oil field development optimization problem where decision parameters of new wells are determined.
In practical settings, the number of drilling rigs is limited and new wells are drilled sequentially, one or a few at a time.
We let $n_{\text{rig}}$ denote the number of rigs (in this work, we set $n_{\text{rig}}$=2), and we let $T_{\text{drill}}$
denote the duration (in days) that it takes to drill and complete a new well and start production/injection.
The initial field development plan is determined by optimizing the (full) development plan over a prior set of geological models,
generated by use of a geostatistical simulation method.
We let $t_i$, $i\ge 0$, denote the time that the $(i+1)$th group of wells is drilled.
The first group of wells is drilled at $t_0=0$.

After drilling and completion of each group of $n_{\text{rig}}$ new wells,
new spatial data (hard data from well locations) and new temporal data (production data)
are collected.
The CLFD model calibration is performed at each $t_i$ ($i\ge 1)$
to update models for consistency with all data.
Here, observed hard data corresponds to the rock property value at (new and existing) well location(s) in the true-model (without any noise added), and is integrated through a geostatistical algorithm.
Synthetic observed production data (for model calibration) is generated by simulating the true-model and adding random noise to the true data.
The true-model here is a realization randomly selected from the prior set.

The optimization is then performed to determine decision parameters
for all future wells (those to be delivered at $t_{i+1}$ and later steps), and controls
of existing wells. Note that only the controls after $t_{i}$ are optimized here.
Optimizing the full development plan at each step is to avoid
a greedy approach where one would only optimize the next $n_{\text{rig}}$ wells
by either freezing or not considering later wells in the project horizon.
A comparison between greedy (well-by-well) optimization and optimization of full development plan is discussed in Example~1 in \citet{shirangi:15b}.
CLFD workflow is shown in Fig.~\ref{fig:schem-clfd}. The optimization and model calibration steps
(together with notations) are discussed next, in turn.

\subsection{CLFD optimization step} \label{opt}
In CLFD, we use $M^i= \{\mathbf{m}_{1}^i,  \ \ldots \ \mathbf{m}_{N_R}^i\}$,
a set of $N_R$ models, to represent the (current) geological knowledge at step $i$.
At each CLFD optimization step, the expected net present value (NPV)
of the project is optimized for a set of representative models selected from the full set.
The problem is formulated as
\begin{equation}\label{eq:opt}
\begin{split}
& \text{maximize} \: \bar{J}(\mathbf{x},M^i_{\text{rep}}), \\
&  {\rm subject~to} \ \  \: \mathbf{x}^{\text{l}}\leq\mathbf{x}\leq \mathbf{x}^{\text{u}}, \quad  g_k(\mathbf{x}) \le 0,\:\: 0\le k \le n_{\text{ie}},
\end{split}
\end{equation}
where $\mathbf{x}$ is the vector of decision parameters,
$M^i_{\text{rep}}= \{\mathbf{m}_{r_1}^i,  \ \ldots \ \mathbf{m}_{r_{n_r}}^i\}$ is a set of $n_r$ representative models ($M^i_{\text{rep}}$ is a subset of $M^i$),
$\mathbf{x}^{\text{l}}$ and $\mathbf{x}^{\text{u}}$ denote the lower and upper bound vectors,
$n_{\text{ie}}$ is the number of inequality constraints,
and $\bar{J}$ is the expected NPV defined as
\begin{equation}\label{eq:npv}
\bar{J}(\mathbf{x},M_{\text{rep}}) = \frac{1}{n_r}\sum_{f=1}^{n_r}J(\mathbf{x},\mathbf{m}_f).
\end{equation}
Computing NPV value for each model $f$ involves a flow simulation run.
In this work, the water-flooding process with two-phase flow is considered where NPV
is computed as
\begin{equation}\label{eq:npv}
\begin{split}
& J(\mathbf{x},\mathbf{m}_f) = \sum_{l=1}^{N_l} \left[\sum_{k=1}^{N_P}(p_{\text{o}}q_{\text{o},k}^l-c_{\text{wp}}q_{\text{w},k}^l)-\sum_{k=1}^{N_I}c_{\text{wi}}q_{\text{wi},k}^l \right] \frac{\Delta t^l}{(1+r_d)^{(t_l/365)}}\\
& -\sum_{w=1}^{N_P+N_I}\frac{C_{\text{well}}}{(1+r_d)^{(t_w/365)}},
\end{split}
\end{equation}
where $C_{\text{well}}$ is the cost of drilling and completing a well,
$t_w$ is the time of delivering well $w$,
$\mathbf{m}_f$ is the vector of reservoir model which contains properties such as permeability values at grid blocks,
$N_l$ is the number of simulation time steps,
$t_l$ is the simulation time (in days),
and $r_d$ is the annual discount rate.
Variables $N_P$ and $N_I$ denote the number
of producers and injectors, respectively, and
$p_{\text{o}}$, $c_{\text{wp}}$ and $c_{\text{wi}}$
indicate the oil price and the cost of
handling produced and injected water (all in $\$/\text{STB}$).
Variables $q_{\text{o},k}^l$ and $q_{\text{w},k}^l$ denote
oil and water production rates for producer $k$
at simulation time step $l$,
$q_{\text{wi},k}^l$ is
the water injection rate of injector $k$ (all in STBD),
and $\Delta{t^l}$ is the size (in days) of simulation time step $l$.

The vector of decision parameters, $\mathbf{x}$, specify
the number of wells, their locations and controls, well type (producer/injector) and the sequence in which new wells are drilled.
The continuous control variables for new and existing wells are treated similarly in the optimization.
For existing wells, only control variables after $t_i$ are included in the optimization at $t_i$.
The location variables are discrete (pseudo-continuous) and treated similar to the continuous variables
but are rounded for input into the reservoir simulation.
The optimization problem in (\ref{eq:opt}) is solved through a parallelized PSO-MADS optimization algorithm \citep{isebor:14b}.
Nonlinear constraint here is the minimum well distance constraint which is handled through the filter method in PSO-MADS.
The PSO-MADS algorithm converges when a maximum number of iterations or a minimum stencil length is reached.
This algorithm is applied within the optimization with sample validation (OSV) procedure, discussed next. 
The initial guess for optimization (through OSV) at $t_i$ corresponds to the optimal solution at $t_{i-1}$.
This is also demonstrated in Fig.~\ref{fig:schem-clfd}.

In CLFD optimization at any $t_i$, the NPV for each model in Eq.~\ref{eq:npv}
is computed from time zero to the end of specified project life.
This is to ensure that we can compare the value of the project at different steps.
Although the NPV is computed from time zero, the simulation runs are only performed for the duration of [$t_i, T$] (as the earlier state variables and decision parameters correspond to the past and do not change).
For restarting the simulations at $t^i$, the state variables for each realization are saved into restart files from the final simulation performed at the model calibration step.
We let $\mathbf{x}^{i}$ denote the optimal solution obtained from the optimization at step~$i$.
Note that other economic measures such as rate of return could be incorporated in the objective function here \citep{shirangi:17b}.

\subsection{Optimization with sample validation and selection of representative models} \label{opt-robust}
We now discuss the optimization with sample validation (OSV) procedure from \citet{shirangi:15b}.
Ideally, the optimization should be performed over the full set of $N_R$ models (here, we set $N_R$ = 50).
Since the computational cost scales linearly with the number of realizations used,
an optimization over the full set will be computationally expensive.
Therefore, we perform the optimization over a set of $n_r << N_R$ representative models (e.g., $n_r=5$).
We will discuss later how these $n_r$ models are selected.

The optimization over the representative subset, improves the initial objective value of $\bar{J}(\mathbf{x}_{\text{init}},M^i_{\text{rep}})$
to an optimal value of $\bar{J}(\mathbf{x}_{\text{opt}},M^i_{\text{rep}})$.
We then apply this optimal solution to the full set of models.
The NPV improvement for the full set of models is therefore
$\bar{J}(\mathbf{x}_{\text{opt}},M^i)$ - $\bar{J}(\mathbf{x}_{\text{init}},M^i)$.
The relative improvement, $RI$, is computed as the ratio of NPV improvement for the full set, divided by that computed for the representative subset, i.e.,
\begin{equation}\label{eq:RI}
RI = \frac{\bar{J}(\mathbf{x}_{\text{opt}},M^i) - \bar{J}(\mathbf{x}_{\text{init}},M^i)}{\bar{J}(\mathbf{x}_{\text{opt}},M^i_{\text{rep}}) - \bar{J}(\mathbf{x}_{\text{init}},M^i_{\text{rep}})}.
\end{equation}

We require $RI\ge \alpha$ to accept $\mathbf{x}_{\text{opt}}$ as the optimal solution at step $i$ ($0 \le \alpha \le 1$)
and we set $\alpha$ = 0.5.
If this criterion is not satisfied, we increase $n_r$
and repeat the optimization, until $RI\ge$ 0.5 is satisfied, or the maximum number of OSV iterations (specified to 4 here) is reached.
In this work, a sequence of $\{5, 9, 16, 25\}$  models is considered in OSV, i.e.,
a PSO-MADS optimization is first performed over $n_r=5$ models with $\mathbf{x}_{\text{opt}}$ as the optimal solution;
if the $RI$ criterion is not satisfied, the optimization is repeated for a (new) set of  $n_r=9$  models selected based on the current solution.
The OSV procedure is outlined as follows.
Note that each OSV iteration, $1\le k \le4$, involves an order of 100 PSO-MADS iterations at Step~3.
\begin{enumerate}
\item At CLFD step $i$, initialize the OSV iteration index to $k=0$, and set $\mathbf{x}_k = \mathbf{x}_{\text{init}} = \mathbf{x}^{i-1}$ where $\mathbf{x}^{i-1}$
is the optimal solution from the previous CLFD optimization step.
Also set the initial value of $n_r$.
\item Select $n_r$ representative models from the full set of $N_R$ realizations in $M^i$ by use of flow simulation results based on $\mathbf{x}_k$.
Note that for $k > 1$, these simulation results are already available from the simulations performed to compute the RI criterion at OSV iteration $k-1$.
\item Solve the optimization problem in (\ref{eq:opt}), using PSO-MADS, with $\mathbf{x}_k$ as the initial guess and the optimal solution denoted by $\mathbf{x}_{k+1}$ ($=\mathbf{x}_{\text{opt}}$).
\item Compute $RI$ in Eq.~\ref{eq:RI} by replacing $\mathbf{x}_{\text{opt}}$ with $\mathbf{x}_{k+1}$ and $\mathbf{x}_{\text{init}}$  with $\mathbf{x}^{i-1}$.
\item If $RI\ge 0.5$, accept the optimal solution and set $\mathbf{x}^i=\mathbf{x}_{k+1}$.\\
If the maximum OSV iterations has reached, set $\mathbf{x}^i = \argmax_{\mathbf{x}_{k}}{\bar{J}(\mathbf{x}_{k},M^i)}$.
Otherwise, set $k=k+1$, increase $n_r$, and go to 2 (note that new representative models will be selected
based on the new $\mathbf{x}_{\text{opt}}$).
\end{enumerate}

The sample validation in OSV is to ensure that the optimal solution obtained for a small representative subset of models
adequately improves the objective function for the full set.
Other criteria such as out-of-sample validation could also be applied and tested here.
For computational results in this work, a sample size of $n_r$ = 9 or 16 typically satisfied $RI\ge 0.5$, while
a smaller size was not usually adequate.

We now discuss the selection of representative models in CLFD.
There are various methods presented for realization selection. Here we use a CDF approach as explained in \citet{shirangi:17}.
In this approach, flow simulation is performed for the full set of $N_R$ models by use of current $\mathbf{x}_k$, and the NPV values are computed.
The cumulative distribution function (CDF) plot is then generated. The models are selected such that two of them
correspond to P10 and P90 (for $n_r \le 10$) or P5 and P95 (for $n_r > 10$), and the rest of $n_r-2$ models correspond to even increments in NPV percentile (e.g., for $n_r=5$, the models correspond to P10, P30, P50, P70, and P90).
Other selection methods such as equally weighted flow-based and permeability-based clustering \citep{shirangi:16b} could be applied and tested here
\footnote{the R code is available here:\\ \url{https://github.com/mehrdad-shirangi/representative_model_selection}}.

\subsection{CLFD model calibration for channelized models} \label{HM-opca}
In CLFD model calibration step,
the $N_R$ realizations are updated with new hard data and production data.
Conditioning to production data is an inverse problem where one or multiple models need to be generated
such that they reproduce the observed flow data and honor prior geological information.
In the context of inverse problem theory \citep{tarantola:05},
generating multiple calibrated (history-matched) models correspond to a sampling of the posterior probability density function (pdf).
A common approach for solving this problem
is the randomized maximum likelihood (RML) \citep{kitanidis:95,shirangi:14}.
With RML, a sample from the posterior pdf is generated by minimizing
an objective function that quantifies the mismatch between observed and predicted data.
This objective function also has a model mismatch term to preserve the prior geological information.
Model calibration can also be accomplished by use of a Kalman filter \citep{ghorbanidehno:15,ghorbanidehno:18}. 

\begin{figure}
\centering
    \includegraphics[width=0.75\textwidth]
    {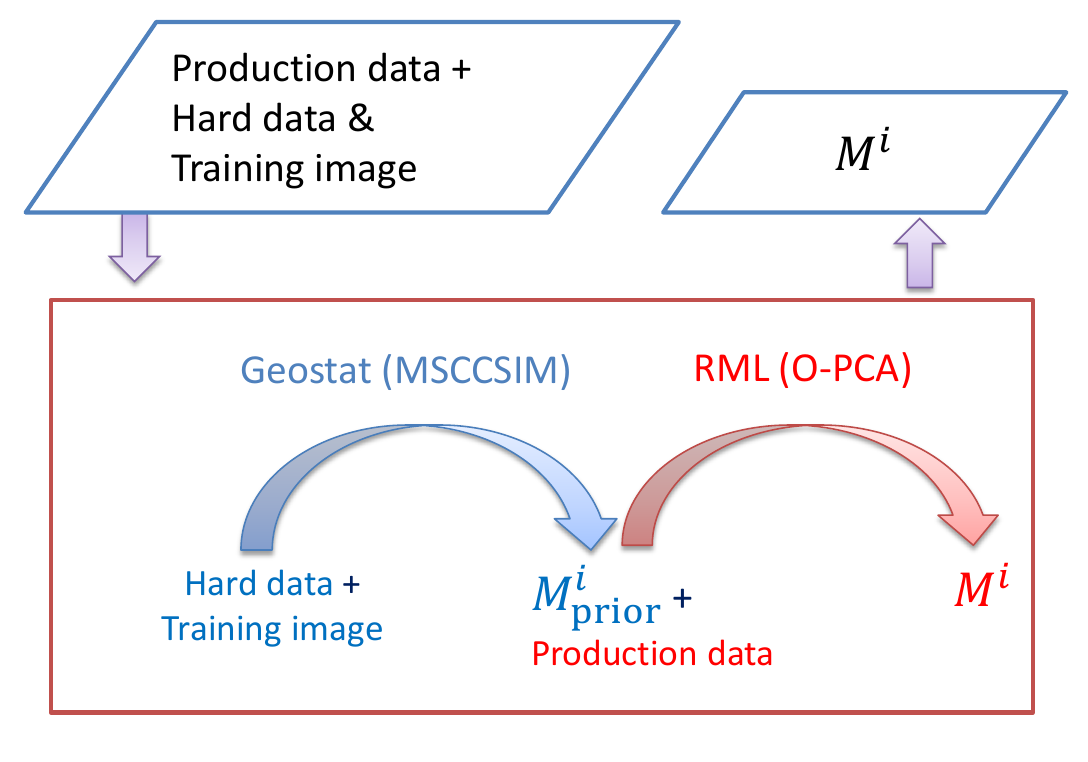}
\caption{Schematic of CLFD model calibration step with multipoint geostatistics.}
\label{fig:schem-hm}
\end{figure}

For model calibration of a channelized reservoir described by MPS,
the optimization-based principal component analysis (O-PCA) parameterization \citep{vo:15,liu:17} is applied here.
The CLFD model calibration for MPS models consists of two steps.
In the first step, new conditional realizations are generated using
a geostatistical simulation approach with hard data from all wells (including the most recent well).
Here we use the multi-scale cross-correlation simulation (MS-CCSIM) geostatistical algorithm \citep{tahmasebi:14} for generating these conditional realizations.
In the second step, the O-PCA-based RML is applied for conditioning
$N_R$ realizations to production data.
The schematic of model calibration procedure is presented in Fig.~\ref{fig:schem-hm}.

We now briefly describe these procedures.
The O-PCA method requires generating $L$ realizations of the permeability field conditioned to hard data at well locations.
We use $L=1000$ in this work.
New realizations must be generated at each CLFD step since new conditioning data become available as we proceed in time.
The MS-CCSIM algorithm is very efficient and can generate these models in a few minutes.
Then the centered matrix of realizations, $\mathbf{X}_{\text{c}}$, is computed,
\begin{equation}  \label{eq:Xc}
\mathbf{X}_{\text{c}} = \begin{bmatrix} \mathbf{m}_1 - \bar{\mathbf{m}} \ldots \: \mathbf{m}_{L}-\bar{\mathbf{m}} \end{bmatrix},
\end{equation}
where $\bar{\mathbf{m}}$ is the mean of the $L$ realizations.
A truncated SVD of $\mathbf{X}_{\text{c}} $
is then computed as $\mathbf{U}_l \mathbf{\Lambda}_l \mathbf{V}_l^T$, where $l<L$.
Given a random $l$-dimensional vector $\boldsymbol{\xi}\sim \mathcal{N} (\mathbf{0},1)$, a new realization
can be generated by solving the following optimization problem:
\begin{equation}  \label{eq:mopca}
\mathbf{m} = \text{argmin}_{\mathbf{z}} \{ \left\Vert \mathbf{U}_l \mathbf{\Lambda}_l \boldsymbol{\xi} + \bar{\mathbf{m}} - \mathbf{z} \right\Vert_2^2 + \gamma R\},
\end{equation}
where $R$ is a regularization term that is specified such that the realization is generated
consistent with the training image (see \citet{vo:14} for details).
For binary models here,
$R$=$\mathbf{z}^T(\mathbf{1}-\mathbf{z})$, where $\mathbf{1}$ is the unity vector of same length as $\mathbf{z}$.
The SNOPT algorithm \citep{gill:05} is applied to solve (\ref{eq:mopca}).

In the O-PCA RML,
a calibrated model is generated by minimizing the following objective function.
\begin{equation}  \label{eq:opc-hm}
S(\boldsymbol{\xi})=\frac{1}{2} (\boldsymbol{\xi}-\boldsymbol{\xi}_{\text{uc}})^T (\boldsymbol{\xi}-\boldsymbol{\xi}_{\text{uc}}) + \frac{1}{2}(\mathbf{g}^p(\boldsymbol{\xi})-\mathbf{d}_{\text{uc}}^{p})^T C_{d,p}^{-1} (\mathbf{g}^p(\boldsymbol{\xi})-\mathbf{d}_{\text{uc}}^{p}),
\end{equation}
where $\boldsymbol{\xi}_{\text{uc}}$ corresponds to a projected MPS realization, i.e.,
$\boldsymbol{\xi}_{\text{uc}} = \mathbf{\Lambda}_l^{-1} \mathbf{U}_l(\mathbf{m}_\text{uc}-\bar{\mathbf{m}})$.
Here $\mathbf{m}_\text{uc}$ is a realization that is unconditioned to production data, but conditioned to hard data ($\mathbf{m}_\text{uc}$
is also generated by MS-CCSIM).
Note that the hard data mismatch term
does not appear since hard data are already honored in the realizations and thus in the O-PCA representation.
Minimization of Eq.~\ref{eq:opc-hm} is also performed by SNOPT \citep{gill:05}. This minimization
is to find the optimal $l$-dimensional $\boldsymbol{\xi}$.
In the computational results of this work, we will set $l = 100$.

Note that for each trial $\boldsymbol{\xi}$, the vector $\mathbf{m}$ is obtained from solving (\ref{eq:mopca}).
Predicted data, $\mathbf{g}^p$, is then generated by performing a reservoir simulation run using this $\mathbf{m}$.
The mismatch objective function in Eq.~\ref{eq:opc-hm} is then computed.
The gradient of $S$ with respect to $\mathbf{m}$ is constructed
through an adjoint solution, using the automatic differentiation framework \citep{bukshtynov:15}.
The gradient is then projected using the chain rule to obtain derivatives
with respect to $\boldsymbol{\xi}$.

\subsection{Parallelized implementation}
The CLFD algorithm can be parallelized in both the optimization and model calibration steps.
In the PSO-MADS algorithm, the objective function for $p$ PSO particles and the $s$ stencil points in MADS algorithm,
can be computed simultaneously, by accessing
$n_{\text{node}}$ compute nodes. The number of stencil points is twice the number of decision parameters, and
$p$ is typically selected to be a comparable number (e.g., $p=s/3$ or $p=s/2$). See \citet{isebor:14a} for more details
on parallelization of PSO-MADS.

In the model calibration step, as the RML runs are independent of one another,
each RML realization is generated on a separate compute node using distributed computing.
All flow simulations are performed using Stanford's Automatic-Differentiation-based General Purpose Research Simulator (AD-GPRS) \citep{younis:11}.
The existing OpenMP-based parallelized version of AD-GPRS \citep{zhou:12m}
allows us to run each simulation on a computational node with 8 cores.
This gives an average speedup of about a factor of 5 for each simulation.
For CLFD results presented in this work,
at each model calibration step, $N_R = 50$ posterior RML models are generated simultaneously using 50 compute nodes, providing a speedup factor
of 250 compared with a sequential application on a single CPU.


\subsection{TruMAP: True-Model-based Assessment of Performance}
As discussed in Introduction, recent work that applied modified versions of CLFD \citep{morosov:16,hanea:17},
reported that application of CLFD, in some cases, has resulted in slight reduction of true-model NPV.
It may seem counter-intuitive that assimilating new data and re-optimizing decision variables over the remaining span of project life,
in cases has lead to a decrease in true-model NPV (which is the ultimate outcome).

It is known in ``decision analysis'' that a decision is distinct from the outcome
 and that a good decision does not necessarily lead to a good outcome \citep{howard:88,bratvold:08}.
In fact, the outcome of a closed-loop optimization procedure can be treated as a random variable.
Therefore, there is a distribution of feasible outcomes at each optimization step.
A sampling of any of these distributions can be obtained by repeating the closed-loop optimization procedure for multiple cases of true-model.
The author refers to this procedure as TruMAP which stands for
``True-Model-based Assessment of Performance''.
The TruMAP procedure is outlined as follows.

\begin{enumerate}
\item Specify $n_{\text{m}}$, the number of true-model cases, and select $n_{\text{m}}$ candidate models, randomly, from the prior set.
\item For each of the $n_{\text{m}}$ true-model cases, perform the closed-loop optimization (CLFD, here) after an optimization over the prior models.
\item Compute all performance measures for each of $n_{\text{m}}$ cases.
\item Present the distribution obtained from these $n_{\text{m}}$ runs for each performance measure.
Following the arguments in \citet{simmons:11}, it is essential not to cherry-pick the results, or eliminate any runs, or reduce $n_{\text{m}}$ here.
\end{enumerate}

We will define performance measures for CLFD and apply the above procedure in the computational results later in this paper.
Note that it is important to
perform the optimization over the prior models for each true-model case as noted in Step~2 of the above procedure.

The use of multiple true-models in CLFD has also been considered by \citet{hanea:17} and \citet{shirangi:15b}.
\citet{hanea:17}, however, did not present a procedure for performance assessment (as TruMAP) and used too few cases.
\citet{barros:16} also applied closed-loop optimization for multiple true-model cases,
though in a different context, i.e., quantifying value of information in CLRM.

The author believes that the TruMAP procedure can help eliminate the unfortunate tendency to
assess an algorithm performance for a few true-model cases and only report the results supporting the significance of the algorithm/hypothesis.
The interested reader is encouraged to refer to \citet{simmons:11} for an excellent discussion and a suggested solution on this.

\section{Computational results and discussion} \label{Ex4}
In this section computational results are presented for a binary channelized system.
In Part~1, detailed application of CLFD is presented.
In Part~2, multiple CLFD runs are performed for performance assessment through TruMAP.

\subsection{Part~1: CLFD for a channelized model} \label{Ex4}
\begin{figure}
\centering
    \includegraphics[width=0.48\textwidth]
    {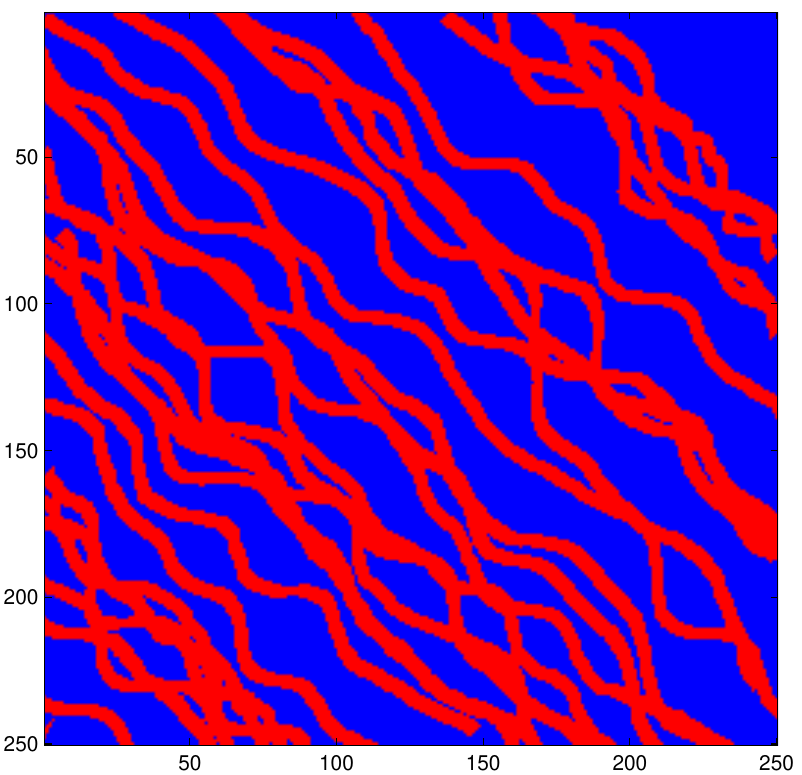}
\caption{Training image used to generate all permeability models (together with any available hard data).}
\label{fig:ti}
\end{figure}

\begin{figure}
\centering
    \includegraphics[width=0.5\textwidth]
    {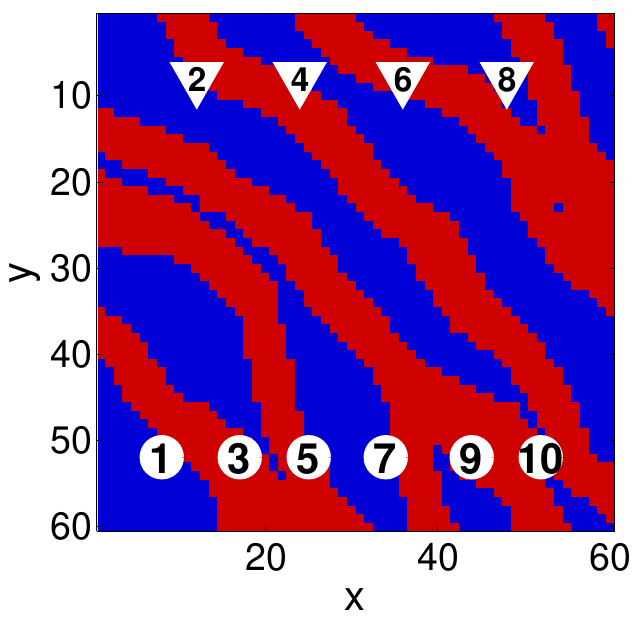}
\caption{True permeability field, with red indicating sand facies (permeability of 500~mD),
and blue indicating shale facies (permeability of 10~mD).
The initial well configuration is also shown with circles denoting producers
and triangles denoting injectors (Part~1).}
\label{fig:ex4truth1}
\end{figure}

\begin{figure}
\centering
   \subfigure[Realization 1]{
    \centering
   \label{fig:ex4real1}
      \includegraphics[width=0.3\textwidth]
      {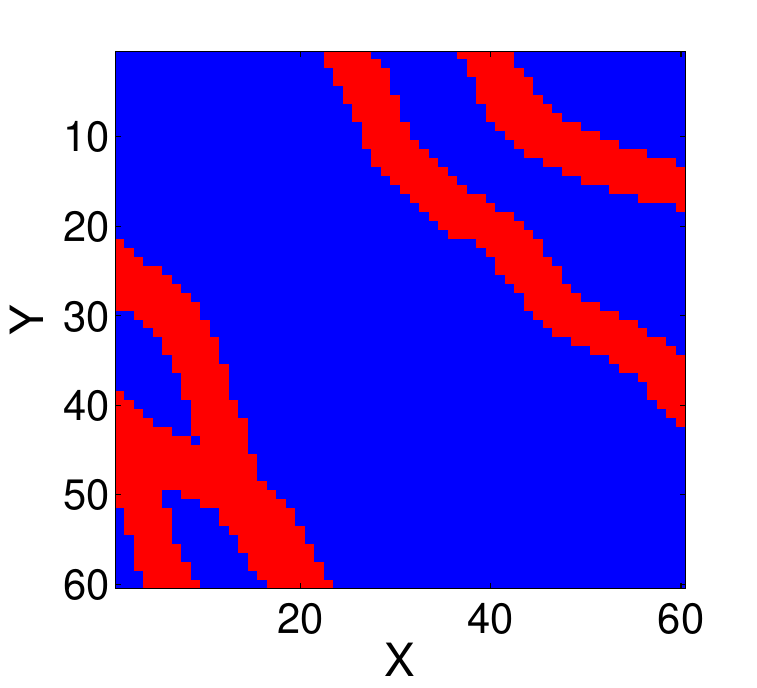}
    }
   \subfigure[Realization 2]{
    \centering
   \label{fig:ex4real2}
      \includegraphics[width=0.3\textwidth]
      {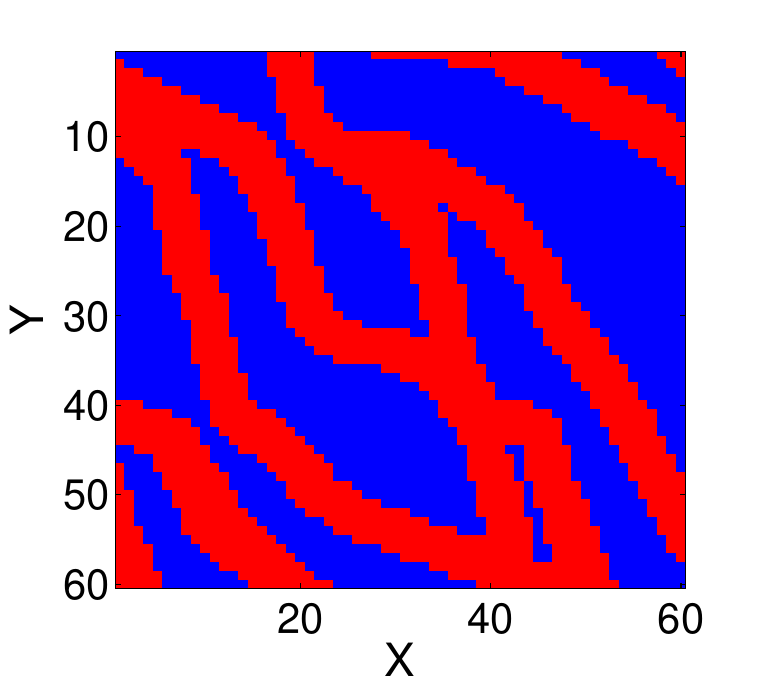}
    }
   \subfigure[Realization 3]{
    \centering
   \label{fig:ex4real3}
      \includegraphics[width=0.3\textwidth]
      {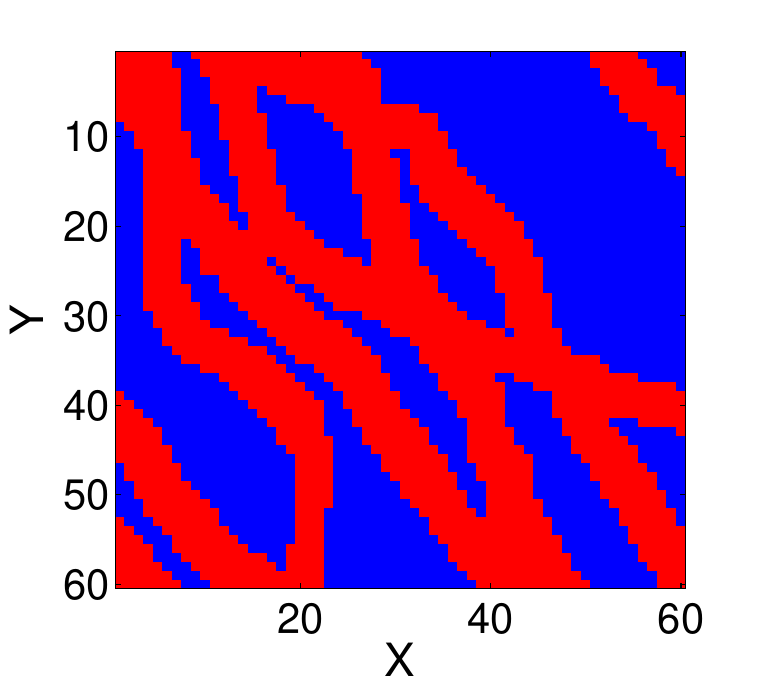}
    }
\caption{Three prior realizations of permeability field, with red indicating sand facies (permeability of 500~mD),
and blue indicating shale facies (permeability of 10~mD).}
\label{fig:ex4prior}
\end{figure}

This example involves a reservoir model described on a two-dimensional uniform grid of dimensions $60\times 60$
with $\Delta{x}=\Delta{y} = 100 \: \text{ft},\: \Delta{z} = 15 \: \text{ft}$.
A binary channelized training image from \citet{vo:14} is taken as the prior geological description (Fig.~\ref{fig:ti}),
from which a set of unconditional realizations are generated using the MS-CCSIM geostatistical algorithm \citep{tahmasebi:14}.
The realizations are not constrained to honor the sand/shale ratio observed in the training image.
The sand permeability is 500~mD, while the shale permeability is 10~mD.
The true permeability field, along with an initial guess for the well locations ($\mathbf{x}^0$),
is shown in Fig.~\ref{fig:ex4truth1}.
Three prior realizations of the permeability field are shown in Fig.~\ref{fig:ex4prior}.
Porosity is assumed to be constant and equal to 0.2.
Initially, the reservoir is at irreducible water saturation which is equal to $S_{wc}$=0.2.

\begin{figure}
\centering
    \includegraphics[width=0.6\textwidth]
    {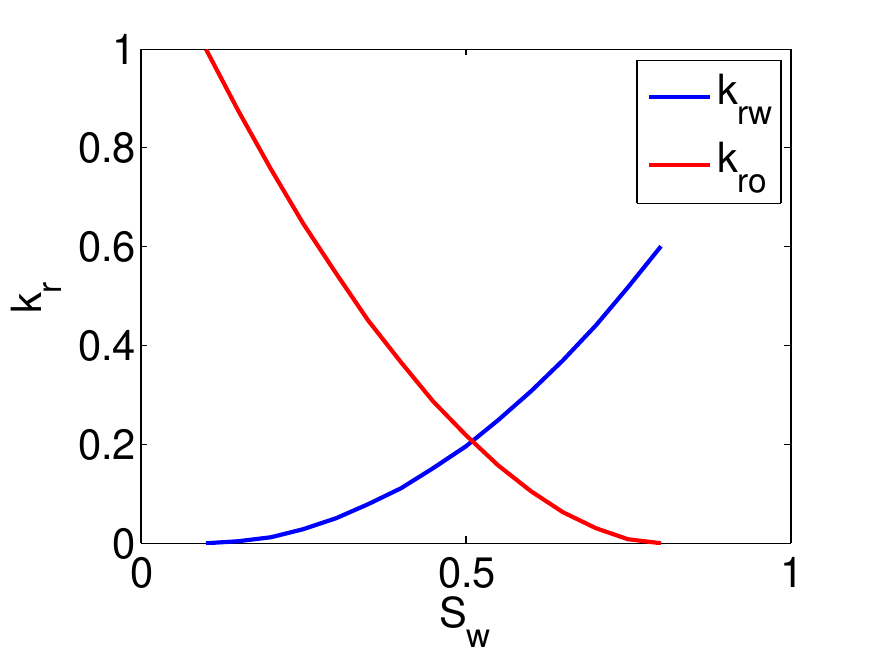}
\caption{Oil and water relative permeability curves. }
\label{fig:relperm}
\end{figure}

Project life is 3000~days, which is divided into seven control steps, with the first five control steps
of length 180~days, and the last two control steps of length 1050 days.
In this example, we assume that two rigs are available ($n_{\text{rig}}$=2) and therefore (a maximum of) two wells can be drilled at each CLFD step.
The optimal number of wells, however, is determined from the optimization.
The drilling time is specified as $T_{\text{drill}}$=180 days. Therefore, the $t_i$ values are given
by $\{0, 180, 360, 540, 720\}$.
The last optimization is performed at 540 days,
which determines the decision parameters corresponding to well type and location of wells 9 \& 10
and (future) BHPs of all wells (wells 1--8).
Discount rate is specified to $r_d=0$.
Simulation and optimization parameters are presented in Table~\ref{Tab:ex1Par},
and relative permeability curves are shown in Fig.~\ref{fig:relperm}.
The simulation input files, permeability realizations, and the training image are available online at
https://purl.stanford.edu/nn001cr1881. 

\begin{table}
\centering \caption{Optimization parameters for all examples}\label{Tab:ex1Par}
\begin{tabular}{lc}
  \hline
Parameter &  Value \\
\noalign{\smallskip} \hline\noalign{\smallskip}
  $C_\text{well}$ & \$10 MM \\
  $r_{o}$ & \$90 STB\\
  $c_{wp}$ & \$10 STB\\
  $c_{wi}$ & \$10 STB\\
  Prod BHP range & $1000-4100$ psi\\
  Inj BHP range  & $4600-7000$ psi\\
\hline\noalign{\smallskip}
\noalign{\smallskip}
\end{tabular}
\end{table}

The number of decision parameters is 80. This corresponds to 10 categorical variables for well types, 20 integer variables
for well locations, and 50 continuous parameters for well settings/controls (BHP).
Not all wells exist at all control steps.
Specifically, in the first control step, only two wells appear in the model.
Similarly, four and six wells exist in the second and third control steps, respectively.
Therefore, the total number of continuous control parameters is 2+4+6+8+10+10+10 or 50.
The initial guess for control parameters correspond to the average of upper bounds and lower bounds.
Well type can take values in $\{-1, 0, 1\}$ which corresponds to \{producer, do not drill, injector\} (for more details please refer to \citet{isebor:14b}).
In each optimization, the PSO-MADS algorithm is applied with 60 PSO particles, and the minimum mesh size for MADS is specified
to be 1\% of the variable range.
In the evaluation of a trial point during optimization,
if the well distance constraint is violated, simulation run is not performed and a zero objective function
value is assigned to that point. This treatment is to avoid frequent failed simulation jobs that happens for trial points with
two wells located too close or at the same grid blocks.
In addition, economic constraint for each producer is enforced by the simulator by specifying a maximum WOR (water to oil ratio).

\begin{figure}
\centering
    \includegraphics[width=0.6\textwidth]
    {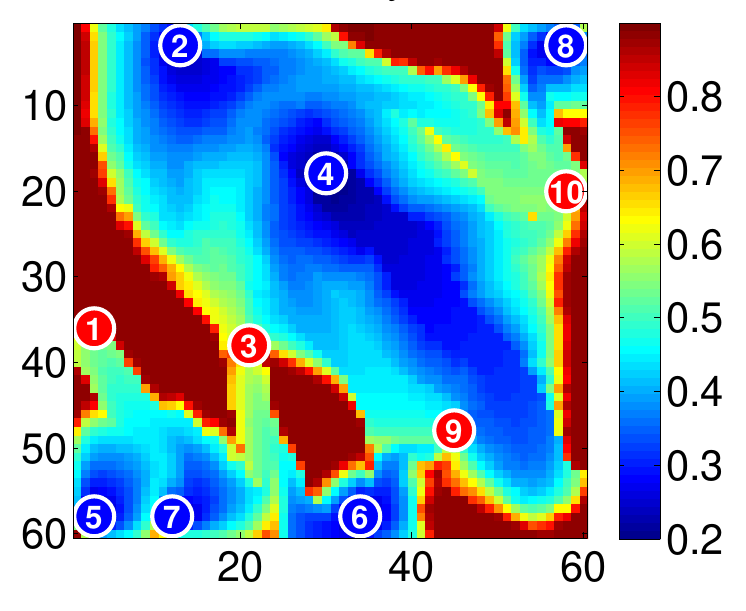}
\caption{
Well configuration from deterministic optimization (using $\mathbf{m}_{\text{true}}$), with red denoting producer, blue denoting injector, and the well numbers indicating the drilling sequence. Background shows final oil saturation.
Note that two wells are drilled at a time. }
\label{fig:ex4swfinal}
\end{figure}

First, deterministic optimization is applied using the true-model (shown in Fig~\ref{fig:ex4truth1}).
The optimal development plan and the final oil saturation map are shown
in Fig.~\ref{fig:ex4swfinal}.
The NPV from the initial guess is \$441 MM while the optimal NPV is \$783 MM.
Next, CLFD using OSV is  applied for this case.

CLFD is applied following the workflow shown in Fig.~\ref{fig:schem-clfd}.
Observed data for CLFD model calibration include production data (oil and water production rates at producers
and water injection rate at injectors) measured at 30-day intervals and hard data from
all existing wells.
Observed production data is generated by adding random noise (measurement error) to the true data,
where the true data corresponds to flow results of the true model when it is simulated with the optimal solution.
Model calibration is performed every 180 days by first generating new geological realizations
conditioned to hard data using the MS-CCSIM algorithm, and then applying O-PCA-based RML using all production data from time 0.

\begin{figure}
\centering
    \includegraphics[width=0.8\textwidth]
    {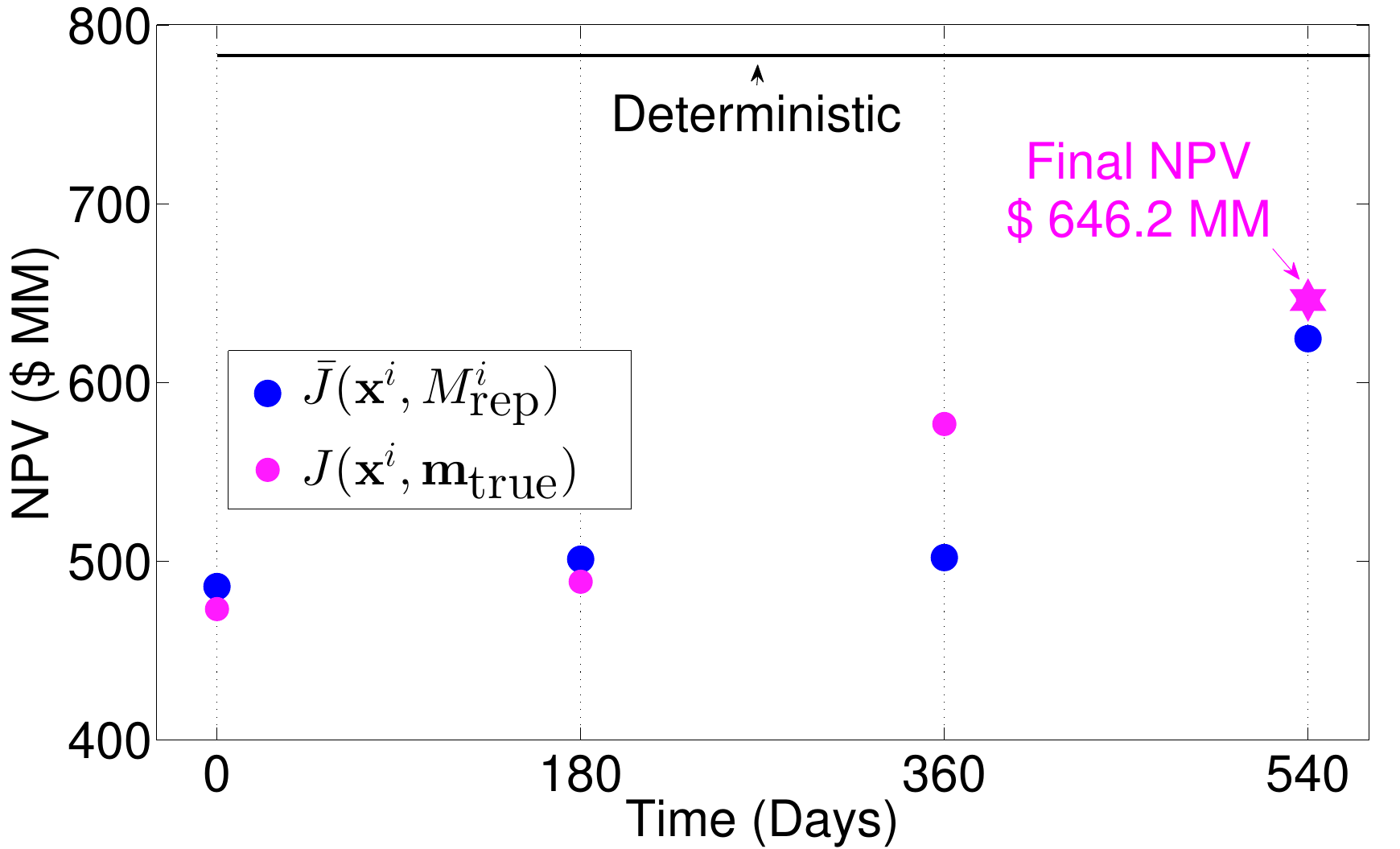}
\caption{Optimal expected NPV, and the corresponding NPV for the true-model, versus CLFD step.
The number of realizations at each CLFD step is determined using OSV. The star shows the final true NPV from CLFD.}
\label{fig:ex4-npvtruth1}
\end{figure}

The progress of the true NPV with CLFD step is shown in Fig.~\ref{fig:ex4-npvtruth1}.
In the evaluation of the true NPV, a nominal strategy is applied where a producer is shut in when the cost of handling
produced water exceeds the revenue from oil.
The NPV from the deterministic optimization is also shown here for comparison.
The final truth-case NPV from CLFD is \$646.2 MM, which is 36.6\% higher than the NPV from optimization (using OSV)
over prior realizations.

\begin{figure}
\centering
    \includegraphics[width=0.8\textwidth]
    {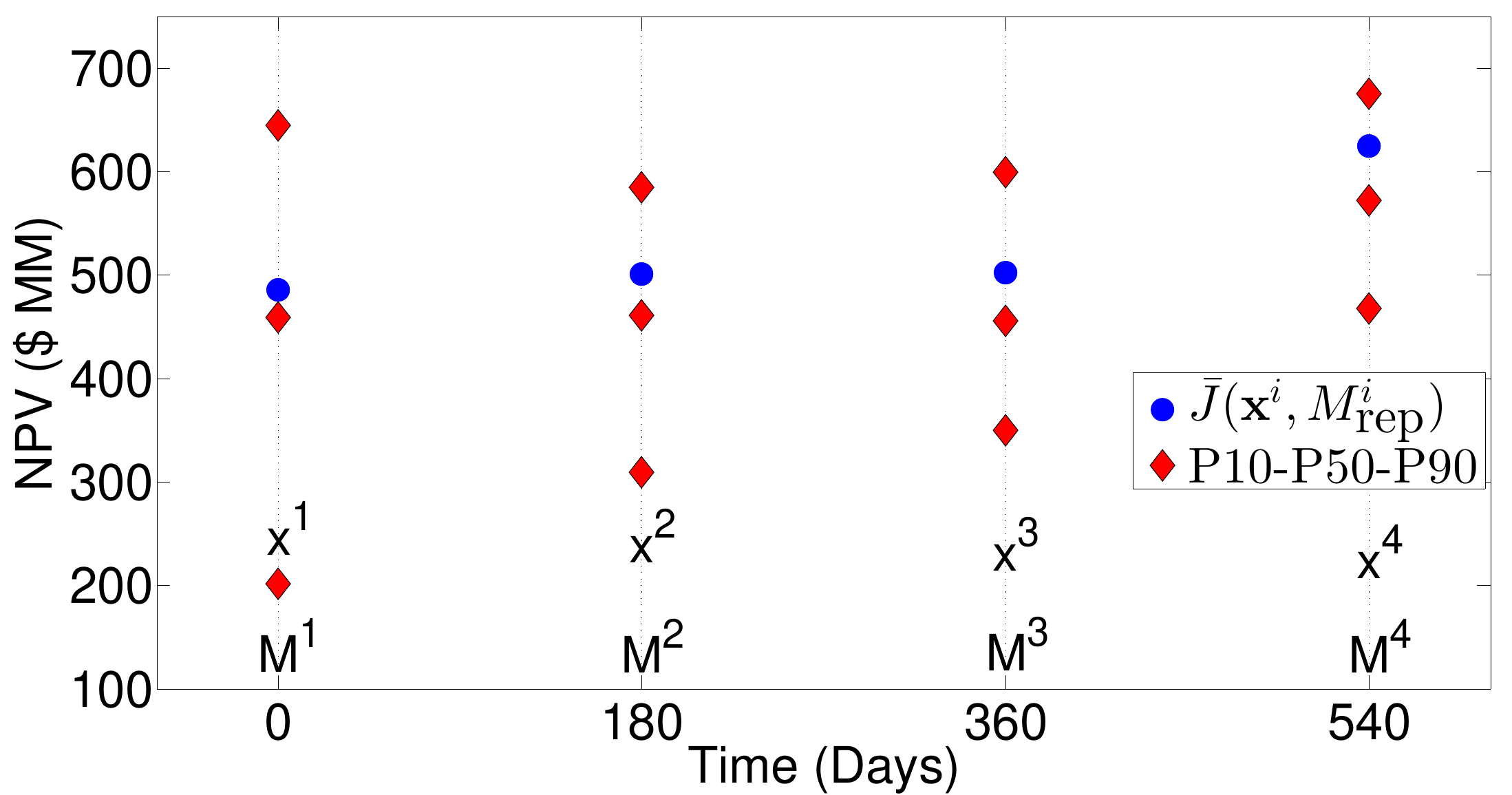}
\caption{P10, P50, P90 NPVs evaluated for the entire set of 50 realizations, along with the expected NPV for the representative subset,
versus CLFD step.}
\label{fig:P10P90}
\end{figure}

Fig.~\ref{fig:P10P90} presents the P10--P50--P90 results for NPV, determined by simulating all 50 realizations and then constructing the cdf, at each CLFD step.
The expected NPV based on the current representative subset (which satisfies validation criterion of $RI \ge 0.5$) is also displayed.
It is evident that the optimal expected NPV for the representative subset falls within the P10--P90 range.
Note that in Fig.~\ref{fig:P10P90}, the percentile values are determined for
the updated models ($M^i$) based on current $\mathbf{x}^i$ at each CLFD step~$i$.
Therefore, the uncertainty range may not necessarily shrink as both the models and the solution are evolving.

\begin{figure}
\centering
   \subfigure[Realization 1 at $t_0$ with $\mathbf{x}^0$]{
    \centering
   \label{fig:real1-x1-Ex4}
      \includegraphics[width=0.3\textwidth]
      {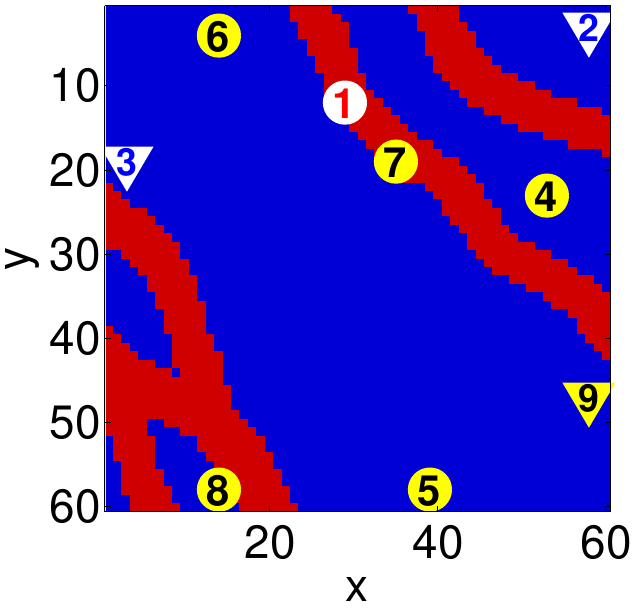}
    }
   \subfigure[Realization 1 at $t_2$ with $\mathbf{x}^2$]{
    \centering
   \label{fig:real1-x3-Ex4}
      \includegraphics[width=0.3\textwidth]
      {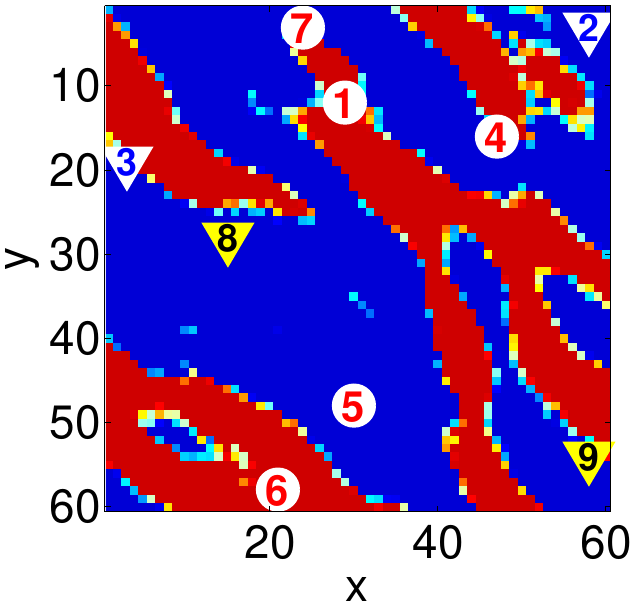}
    }
   \subfigure[Realization 1 at $t_3$ with $\mathbf{x}^3$]{
    \centering
   \label{fig:real1-x4-Ex4}
      \includegraphics[width=0.3\textwidth]
      {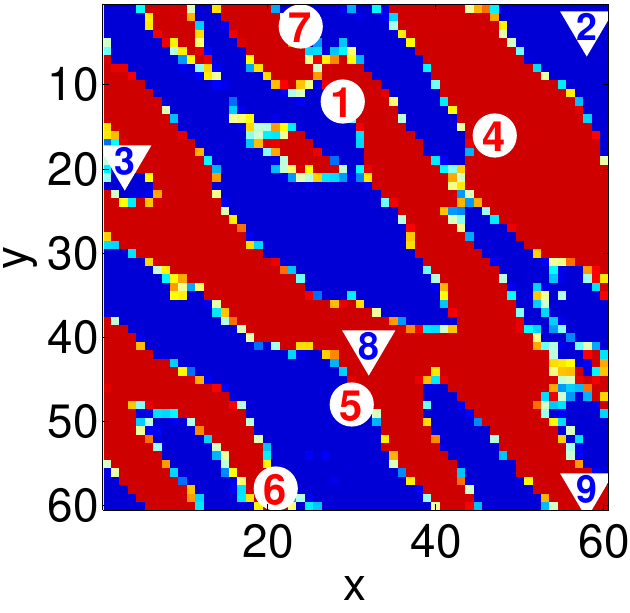}
    }
   \subfigure[Realization 2 at $t_0$ with $\mathbf{x}^0$]{
    \centering
   \label{fig:real2-x1-Ex4}
      \includegraphics[width=0.3\textwidth]
      {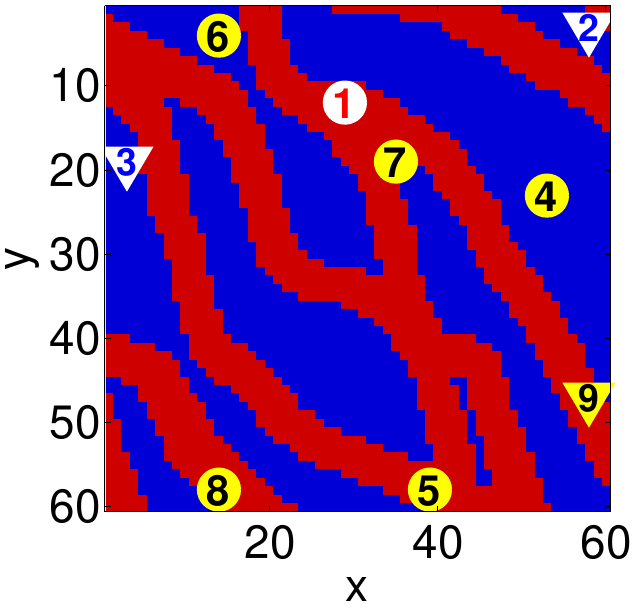}
    }
   \subfigure[Realization 2 at $t_2$ with $\mathbf{x}^2$]{
    \centering
   \label{fig:real2-x3-Ex4}
      \includegraphics[width=0.3\textwidth]
      {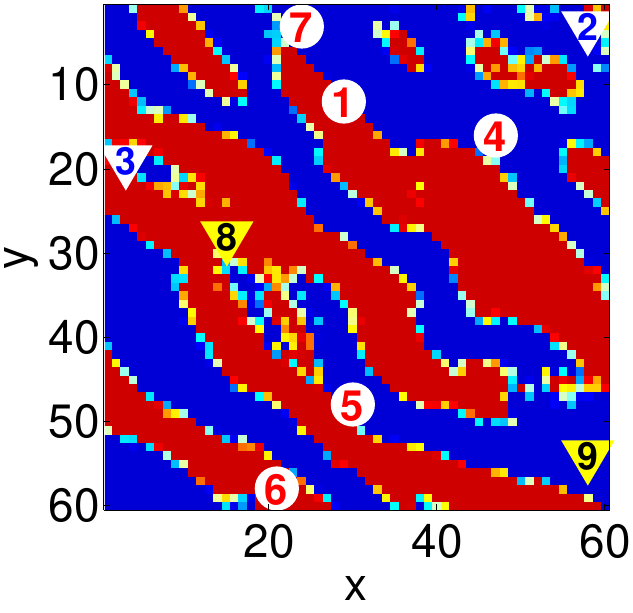}
    }
   \subfigure[Realization 2 at $t_3$ with $\mathbf{x}^3$]{
    \centering
   \label{fig:real2-x4-Ex4}
      \includegraphics[width=0.3\textwidth]
      {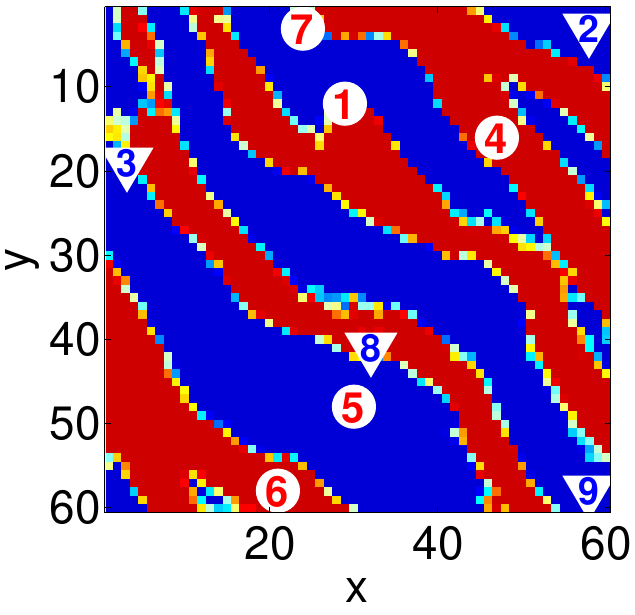}
    }
\caption{Evolution of two RML realizations for different CLFD steps, with red indicating sand facies (permeability of 500~mD),
and blue indicating shale facies (permeability of 10~mD).
Current optimal well configuration and drilling sequence is also depicted. Solid white circles
and triangles denote producers and injectors (drilled or in the process of being drilled),
and yellow circles and triangles denote planned producers and injectors. Numbers indicate the drilling sequence (Example~1).
}
\label{fig:evol-OSV-OPCA}
\end{figure}

Fig.~\ref{fig:evol-OSV-OPCA} shows the evolution of the well configuration and the geological model for two realizations.
The well scenario involves three injectors at $t_1$, but four injectors at later times.
Each realization continues to show differences though the CLFD steps due to conditioning to new hard and production data.

\begin{figure}
\centering
   \subfigure[]{
    \centering
   \label{fig:ex4real1}
      \includegraphics[width=0.45\textwidth]
      {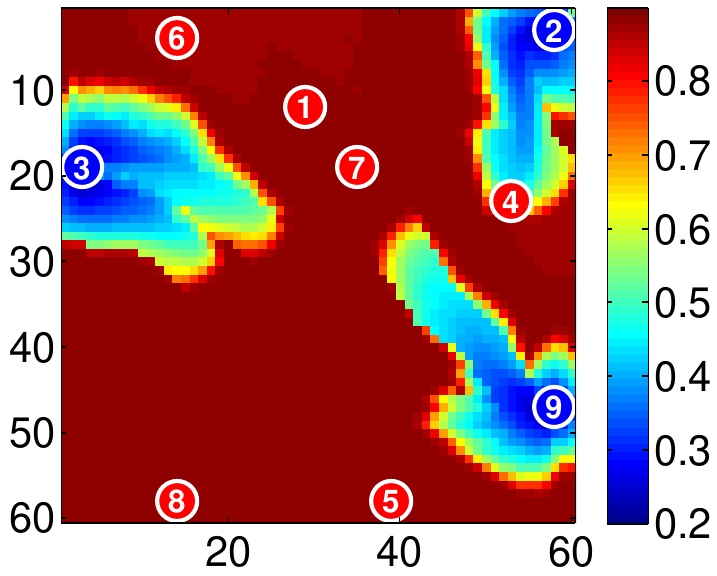}
    }
   \subfigure[]{
    \centering
   \label{fig:ex4real2}
      \includegraphics[width=0.45\textwidth]
      {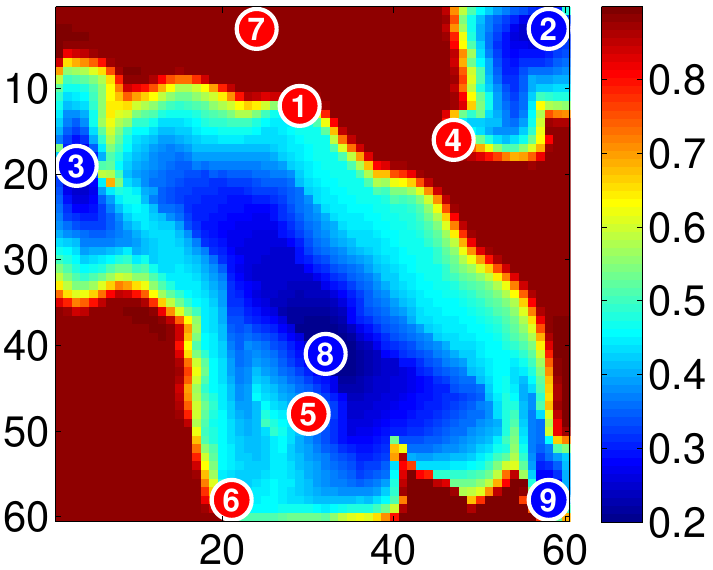}
    }
\caption{Final oil saturation from (a) optimization over prior models by OSV, and (b) final CLFD solution.
Corresponding well configuration is also shown, with red denoting producer, blue denoting injector, and the well numbers indicating the drilling sequence. }
\label{fig:oilsat}
\end{figure}

\begin{figure}
\centering
   \subfigure[]{
    \centering
   \label{fig:npv}
      \includegraphics[width=0.47\textwidth]
      {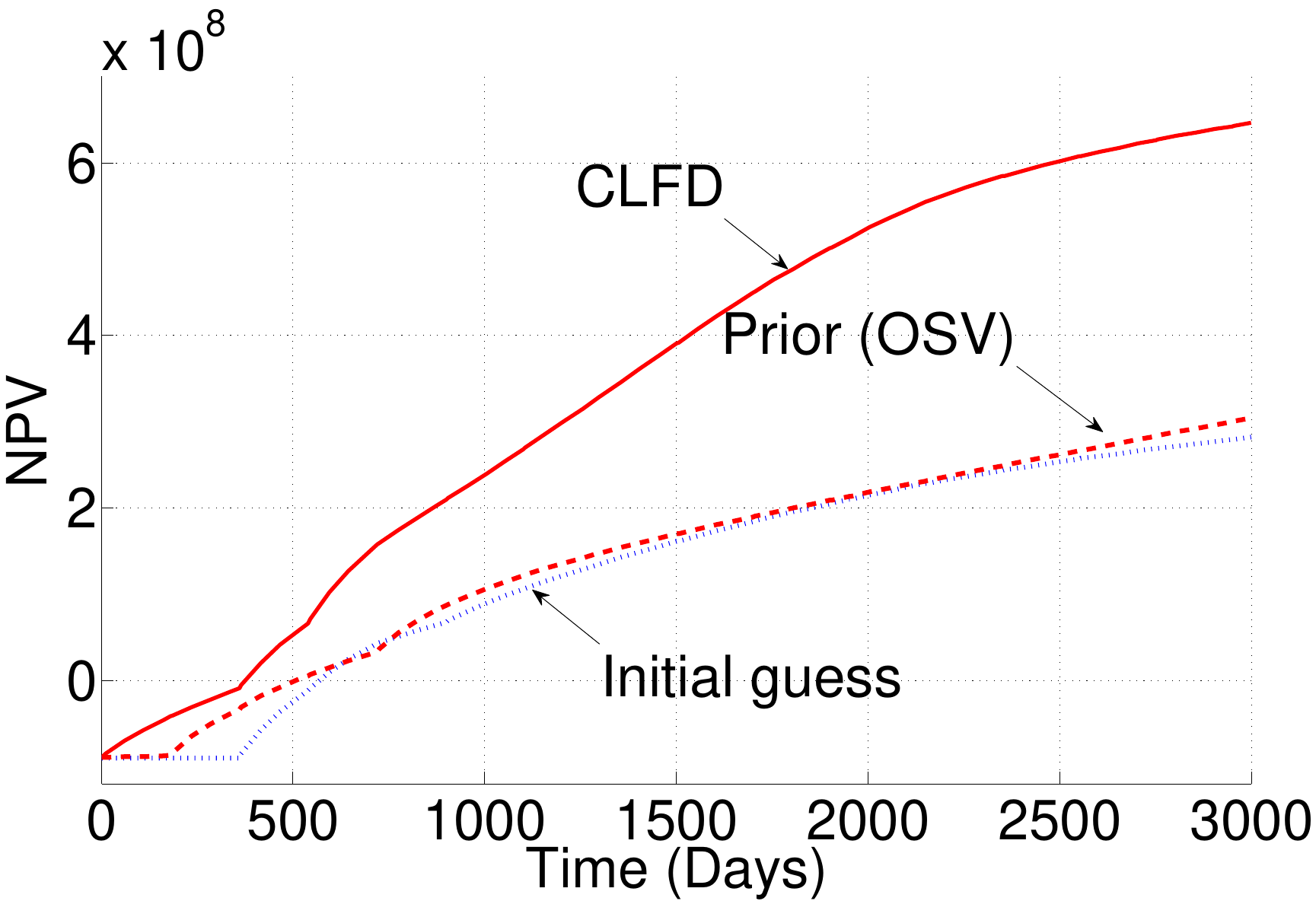}
    }
   \subfigure[]{
    \centering
   \label{fig:cumoilsub}
      \includegraphics[width=0.47\textwidth]
      {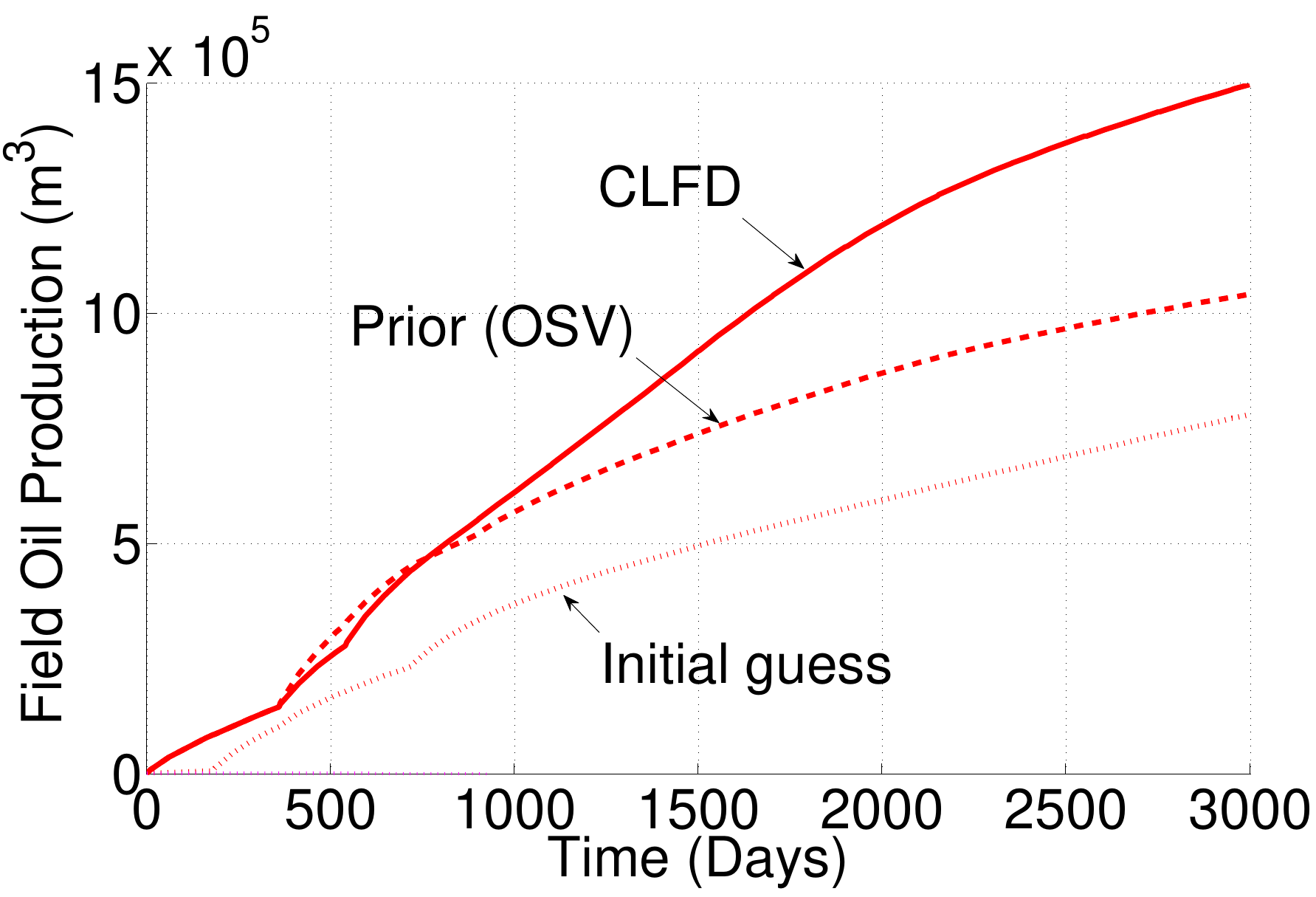}
    }
\caption{(a) NPV trajectory, and (b) cumulative oil production, from different approaches.}
\label{fig:cumOil}
\end{figure}

It may be worth comparing the final oil saturation from final CLFD solution and from the solution obtained by optimization
over prior models through OSV (Fig.~\ref{fig:oilsat}). It is evident that the solution from CLFD corresponds to improved sweep of oil in this case.
Further, the CLFD solution also corresponds to both higher oil production and NPV as shown in Fig.~\ref{fig:cumOil}.
It appears that the optimal solution from OSV only slightly improves the final NPV compared to the initial guess,
whereas the CLFD solution displays a significant improvement.

\begin{figure}
\centering
   \subfigure[ $\bar{\mathbf{m}}_{\text{prior}}$, 180 days]{
    \centering
   \label{fig:mprior-t2}
      \includegraphics[width=0.21\textwidth]
      {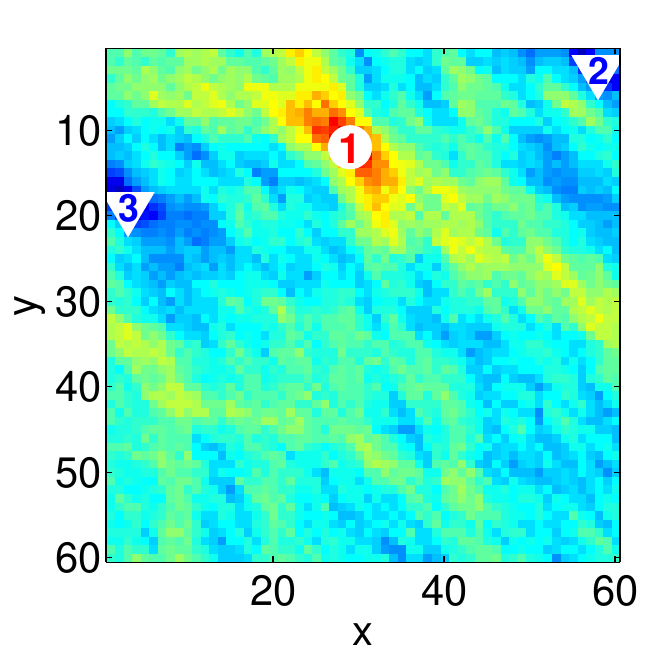}
    }
   \subfigure[ $\bar{\mathbf{m}}_{\text{prior}}$, 360 days]{
    \centering
   \label{fig:mprior-t3}
      \includegraphics[width=0.21\textwidth]
      {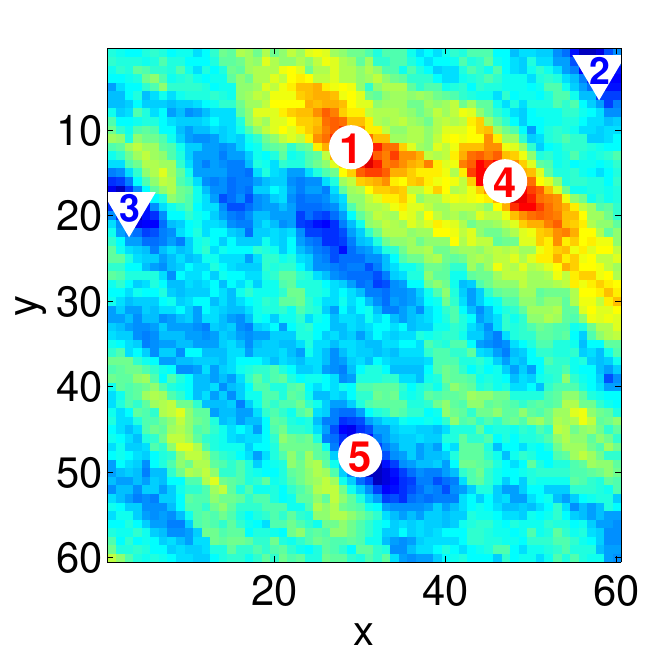}
    }
   \subfigure[ $\bar{\mathbf{m}}_{\text{prior}}$, 540 days]{
    \centering
   \label{fig:mprior-t4}
      \includegraphics[width=0.21\textwidth]
      {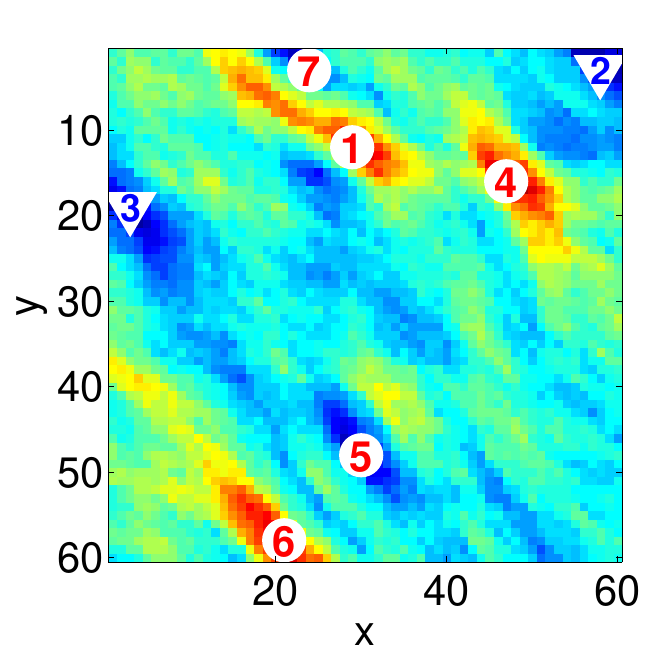}
    }
   \subfigure[ $\bar{\mathbf{m}}_{\text{prior}}$, 720 days]{
    \centering
   \label{fig:mprior-t5}
      \includegraphics[width=0.24\textwidth,height=1.26in]
      {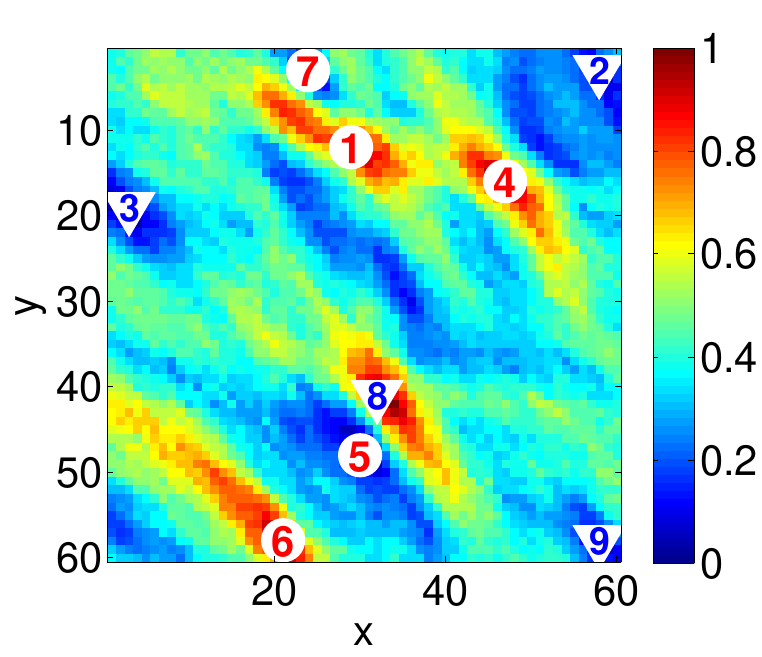}
    }
   \subfigure[ $\bar{\mathbf{m}}_{\text{post}}$, 180 days]{
    \centering
   \label{fig:mpost-t2}
      \includegraphics[width=0.21\textwidth]
      {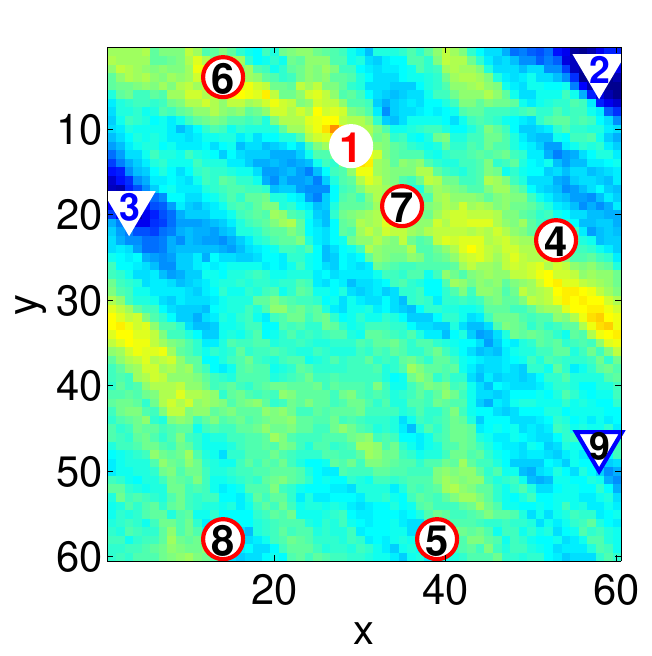}
    }
   \subfigure[ $\bar{\mathbf{m}}_{\text{post}}$, 360 days]{
    \centering
   \label{fig:mpost-t3}
      \includegraphics[width=0.21\textwidth]
      {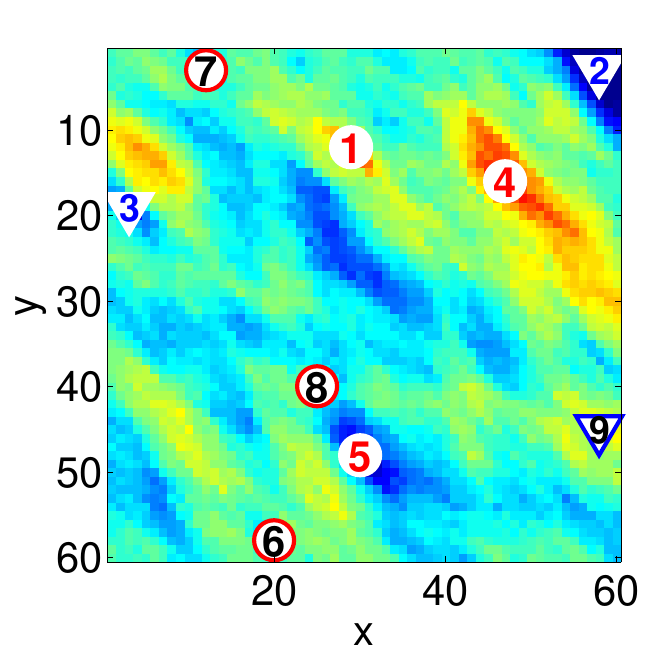}
    }
    \subfigure[ $\bar{\mathbf{m}}_{\text{post}}$, 540 days]{
    \centering
   \label{fig:mpost-t4}
      \includegraphics[width=0.21\textwidth]
      {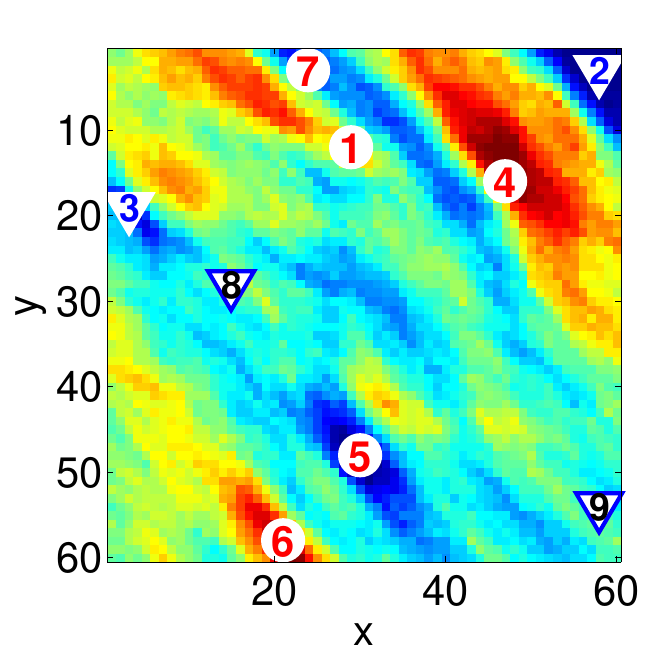}
    }
    \subfigure[ $\bar{\mathbf{m}}_{\text{post}}$, 720 days]{
    \centering
   \label{fig:mpost-t5}
      \includegraphics[width=0.24\textwidth,height=1.26in]
      {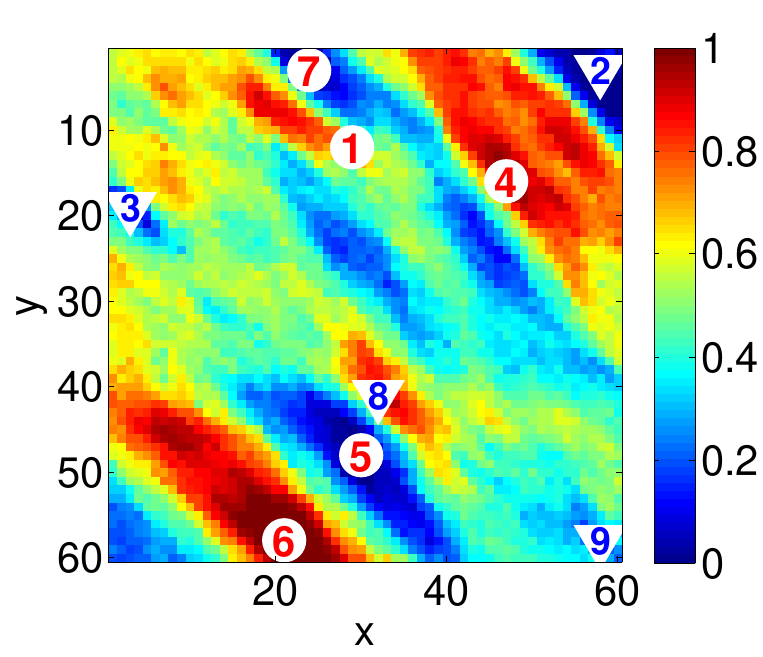}
    }
\caption{Evolution of mean of ($N_R = 50$) prior realizations (conditioned to hard data) and mean of ($N_R = 50$) posterior realizations
of facies distribution, for different CLFD steps.
Current optimal well configuration and drilling sequence is also depicted.
White circles and triangles denote producers and injectors, respectively. Wells with colored (red or blue) numbers are drilled,
while outlined red circles and blue triangles denote planned producers and injectors.
For the prior model (a-d) only the drilled wells are shown.
Numbers indicate the drilling sequence.}
\label{fig:ex4-evolmean}
\end{figure}

As discussed earlier,
the CLFD model calibration step involves integrating production data from all previous wells (except the most recent wells)
together with hard data from all wells including the most recent wells.
Integrating both hard data and production data is required to achieve optimal CLFD performance.
At each CLFD model calibration step, new realizations conditioned to hard data are generated.
We let $\bar{\mathbf{m}}_{\text{prior}}$ designate the mean
of these $N_R = 50$ realizations (conditioned to all available hard data) at each update step.
Figs.~\ref{fig:mprior-t2}-\ref{fig:mprior-t5}
show the evolution of the prior mean with CLFD step.
We also compute the mean of $N_R = 50$ posterior realizations (conditioned to both production and hard data).
These are shown in Figs.~\ref{fig:mpost-t2}-\ref{fig:mpost-t5}.
The optimal development plan at each CLFD step is also shown.
It is evident that the mean of the posterior realizations (conditioned to production data)
more closely resembles the true-model (Fig.~\ref{fig:ex4truth1}).
This indicates the importance of integrating production data to reduce the uncertainty in the geological description.

In terms of computational cost, this CLFD application involved a total number of 356,500 simulation runs.
As computations were performed on a cluster with access to maximum 400 compute cores, the parallelized computations is equivalent
to 1,151 sequential simulation runs (calls to cluster) where each batch of simulations takes about 104 seconds on average.
The total computational time is about 32 hours. Note that in practice, each CLFD step is performed
at a different time and thus the computational cost will be a fraction of the total cost here.

\subsection{Part 2: TruMAP for Statistical assessment of CLFD performance} \label{results:multiruns}
In this section, our objective is to assess the CLFD and the OSV performance through application with multiple truth-case models.
Following the TruMAP procedure, we specify $n_{\text{m}} =25$ and consider 24 additional $\mathbf{m}_{\text{true}}$ models.
We repeat the CLFD application for each of these true-models
so as to assess the performance in terms of various measures as follows.

\subsubsection{Optimization with sample validation results} \label{results:multiruns}
We first discuss the OSV results. We only consider the optimization over prior models ($i=0$) here.
The OSV procedure systematically determines adequate number of representative models in optimization over the full set of realizations (which is 50 here).
The results for five of the $n_{\text{m}}$=25 cases are shown in Table~\ref{tab:osv-pr} (the first row corresponds to the true-model discussed earlier).
The average and median $n_r$ in this case are 12.5 and 9, respectively (Table~\ref{tab:osv-pr}).
We conclude that for this problem,
$n_r=5$ is typically not sufficient, whereas optimization over $n_r=9$ or $16$ realizations is more likely to satisfy the OSV validation criterion.
We also observe that the average true-model NPV with OSV is 15.6\% greater than that with $n_r=5$ prior models.

\begin{table}
\centering \caption{
NPV (\$ MM) values from optimization over $n_r = 5$ prior realizations and by use of OSV (where $n_r$ is increased
to satisfy $RI>0.5$) for prior realizations, for five different true models.
}\label{tab:osv-pr}
\begin{tabular}{lcccc}
  \hline
True model & $J(\mathbf{m}_{\text{true}},\mathbf{x}_{\text{opt}})$ with $n_r=5$ & $J(\mathbf{m}_{\text{true}},\mathbf{x}^{i=0})$ from OSV (RI$\ge$ 0.5) & final $n_r$ \\
\noalign{\smallskip} \hline\noalign{\smallskip}
1	& 474	& 473	& 9 \\
6	& 462	& 376	& 16 \\
11	& 486	& 461	& 9 \\
16	& 402	& 553	& 9 \\
21	& 318	& 323	& 9 \\
Avg (1-25)   & 339  & 392  & 12.5\\
Median (1-25) & -- & -- &   9 \\
\hline\noalign{\smallskip}
\noalign{\smallskip}
\end{tabular}
\end{table}



\subsubsection{Assessment of single-step CLFD performance} \label{results:multiruns}
We now assess the performance of CLFD in terms of improvement of true-model NPV from a single CLFD step.
The NPV values for 5 models at each CLFD step are compiled in Table~\ref{npv-clfd}.
The average true NPV from 25 runs (also shown in Table~\ref{npv-clfd}) increases monotonically with CLFD step.
The average true NPV at CLFD step~1 is 9.9\% greater than optimization over prior ($i=0$),
and the increase from step~1 to 2 and from step~2 to 3 is 7.4\% and 10.1\%, respectively.
At each CLFD step~$i$ ($i\ge 1$), we define the NPV improvement for true-model $j$ as
\begin{equation}  \label{eq:dJi}
dJ_{\text{true},j}^{i} =  \frac{J(\mathbf{m}_{\text{true},j},\mathbf{x}^i_j) - J(\mathbf{m}_{\text{true},j},\mathbf{x}^{i-1}_j)}{J(\mathbf{m}_{\text{true},j},\mathbf{x}^{i-1}_j)},
\end{equation}
where the numerator is the NPV improvement for true-model $j$ as a result of performing one step of CLFD (step~$i$).
The mean and percentile values of $dJ_{\text{true},j}^{i}$ are compiled in Table~\ref{tab:dJ}.
The mean improvement varies between 10--14\%, and the P10 monotonically increase with CLFD step (Table~\ref{tab:dJ}).

\begin{figure}
\centering
    \includegraphics[width=0.8\textwidth]
    {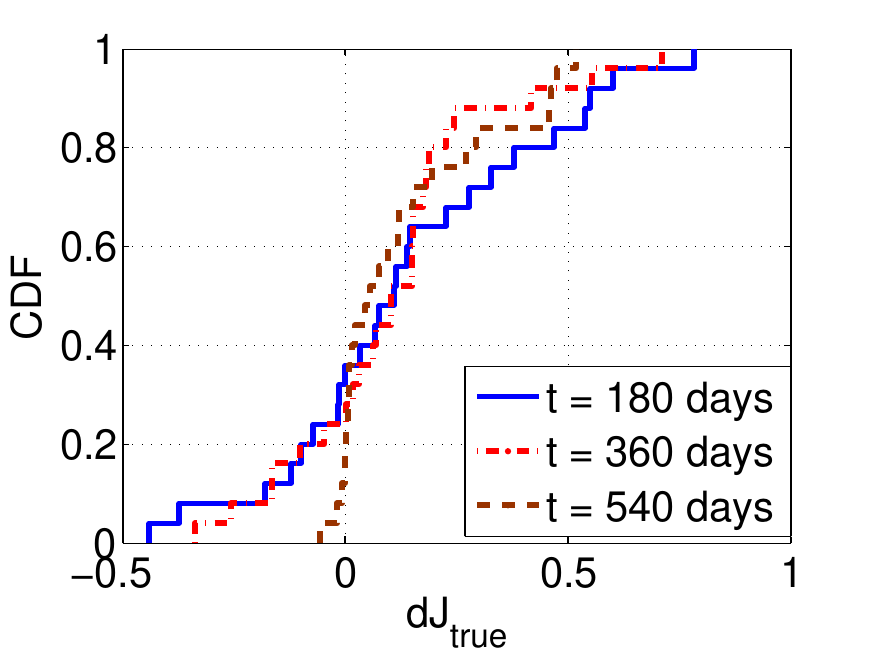}
\caption{Cumulative distribution function (CDF) for 25 $dJ_{\text{true},j}^{i}$ values (relative improvement in true-model NPV from one-step of CLFD run)
for $i$=1 ($t$=180~days), $i$=2 ($t$=360~days), and $i$=3 ($t$=540~days).}
\label{fig:cdf-iStep2}
\end{figure}

The probability of increasing true-model NPV through a single-step of CLFD is simply the fraction of
positive $dJ_{\text{true},j}^{i}$ values (last column of Table~\ref{tab:dJ}).
CLFD improved the true-model NPV for 64\% of cases at $i=1$,
76\% of cases at $i=2$, and 80\% of cases at step $i=3$.
This performance improvement can also be observed visually in Fig.~\ref{fig:cdf-iStep2}.
Therefore, we observe that although application of a single-step CLFD may even reduce the true NPV of the project in some cases,
the chance of increasing true NPV increases at later steps.

\begin{table}
\centering \caption{
NPV values ($J(\mathbf{m}_{\text{true},j})$ in \$ MM) from OSV over prior models
and at each CLFD step.
Average values are also provided.}\label{npv-clfd}
\begin{tabular}{lcccc}
  \hline
True model & OSV-prior ($i$=0) & $t$=180 ($i$=1) & $t$=360 ($i$=2) & $t$=540 ($i$=3)\\
\noalign{\smallskip} \hline\noalign{\smallskip}
1	&  473	&  488	&  577	&  646 \\
6	&  376	&  376	&  314	&  406 \\
11	&  461	&  289	&  449	&  570 \\
16	&  553	&  498	&  529	&  532 \\
21	&  323	&  518	&  596	&  602 \\
Avg (1-25)	&  392	&  431	&  463	&  510 \\
\hline\noalign{\smallskip}
\noalign{\smallskip}
\end{tabular}
\end{table}

\begin{table}
\centering \caption{
Statistics of $dJ_{\text{true}}^{i}$ (true NPV improvement after one step of CLFD) for 25 different true models}\label{tab:dJ}
\begin{tabular}{lccccc}
  \hline
CLFD step (time) & mean  & P10 & P50 & P90 & fraction positive\\
\noalign{\smallskip} \hline\noalign{\smallskip}
$i=1$ (180~days)	& 0.14	& -0.18	& 0.11	& 0.55 & 0.64\\
$i=2$ (360~days)	& 0.10	& -0.16	& 0.10	& 0.42 & 0.76\\
$i=3$ (540~days)	& 0.13	& -0.01	& 0.06	& 0.46 & 0.80\\
\hline\noalign{\smallskip}
\noalign{\smallskip}
\end{tabular}
\end{table}

At $i=1$ ($t=180$ days), the uncertainty reduction in model description
through integration of spatial data (hard data from four wells) and temporal data (production data from two wells)
is not as significant as later steps.
As a result, there is 36\% chance of lowering true-model NPV through a single CLFD step.
At later steps, however, more wells are drilled and more data results in improved model description.
Therefore, there is a higher chance that the optimization over the updated models would result in a higher project outcome (true NPV).


\subsubsection{Assessment of multi-steps CLFD performance} \label{results:multiruns}
To further assess the overall net impact of performing multiple steps of CLFD, we define the relative improvement with respect
to optimization over prior as
\begin{equation}  \label{eq:dJp}
dJ_{\text{rel-pr},j}^{i} =  \frac{J(\mathbf{m}_{\text{true},j},\mathbf{x}^i_j) - J(\mathbf{m}_{\text{true},j},\mathbf{x}^{i=0}_j)}{J(\mathbf{m}_{\text{true},j},\mathbf{x}^{i=0}_j)},
\end{equation}
where rel-pr stands for ``relative to prior''.
For $i$=1, $dJ$ in Eq.~\ref{eq:dJp} provides the same quantity as that in Eq.~\ref{eq:dJi}.
However, at later steps ($i$=2 and $i$=3), these two quantities are different.

\begin{figure}
\centering
    \includegraphics[width=0.8\textwidth]
    {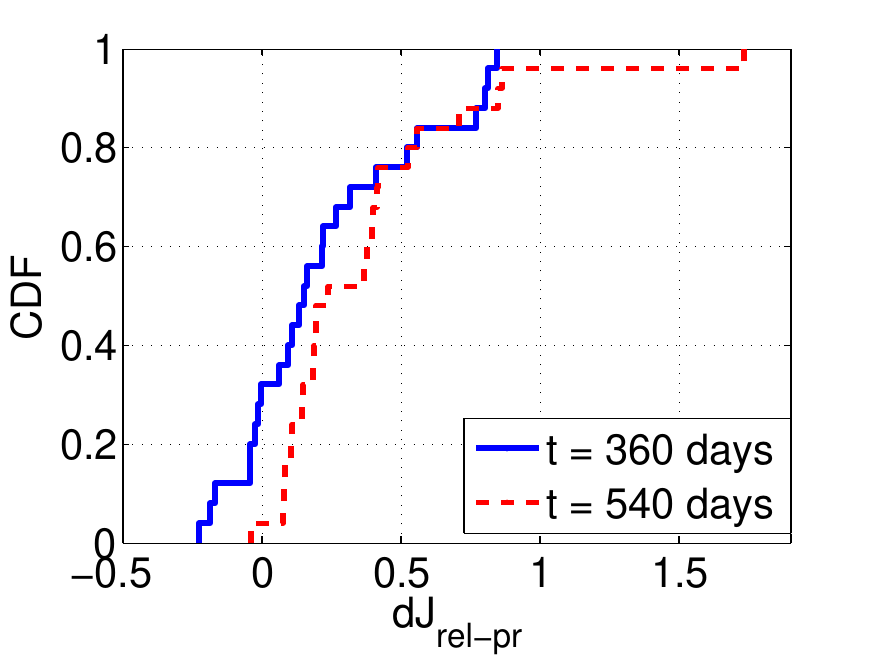}
\caption{Cumulative distribution function (CDF) for 25 $dJ_{\text{rel-pr}}^{i}$ values (improvement in true-model NPV relative to optimization over prior models).}
\label{fig:cdf-dJpr}
\end{figure}

Fig.~\ref{fig:cdf-dJpr} shows the CDF plot for $dJ_{\text{rel-pr}}^{i}$ values, and Table~\ref{tab:dJpr} provides the numeric percentiles.
These results indicate that application of two CLFD steps,
has improved the NPV of 68\% of the models, with an average improvement of 23\%.
The full CLFD run (three CLFD steps in this case), improved the true NPV in 96\% of cases
with an average improvement of 37\%.

These results together with those in the previous section (and Table~\ref{tab:dJ}),
indicate that the probability of increasing true model NPV increases with CLFD step.
In other words, there is a higher chance of improving true NPV through multi-step CLFD,
compared with a single step of CLFD.
In this case, we observe that
application of one, two, and three consecutive steps of CLFD improved the NPV in 64\%, 68\%,
and 96\% of case, respectively.
These results infer about an implicit correction mechanism in CLFD that decisions
about planned wells (those to be drilled at later steps in the time horizon) and future controls are ultimately improved.

\begin{table}
\centering \caption{
Statistics of $dJ_{\text{rel-pr}}^{i}$ (true-model NPV improvement at each CLFD step $i$ relative to optimization over prior models) for 25 different true models}\label{tab:dJpr}
\begin{tabular}{lccccc}
  \hline
CLFD step (time) & mean  & P10 & P50 & P90 & fraction positive \\
\noalign{\smallskip} \hline\noalign{\smallskip}
$i=1$ (180~days)	& 0.14	& -0.18	& 0.11	& 0.55 & 0.64\\
$i=2$ (360~days)	& 0.23	& -0.17	& 0.15	& 0.80 & 0.68\\
$i=3$ (540~days)	& 0.37	& 0.08	& 0.24	& 0.85 & 0.96\\
\hline\noalign{\smallskip}
\noalign{\smallskip}
\end{tabular}
\end{table}

\begin{figure}
\centering
    \includegraphics[width=0.8\textwidth]
    {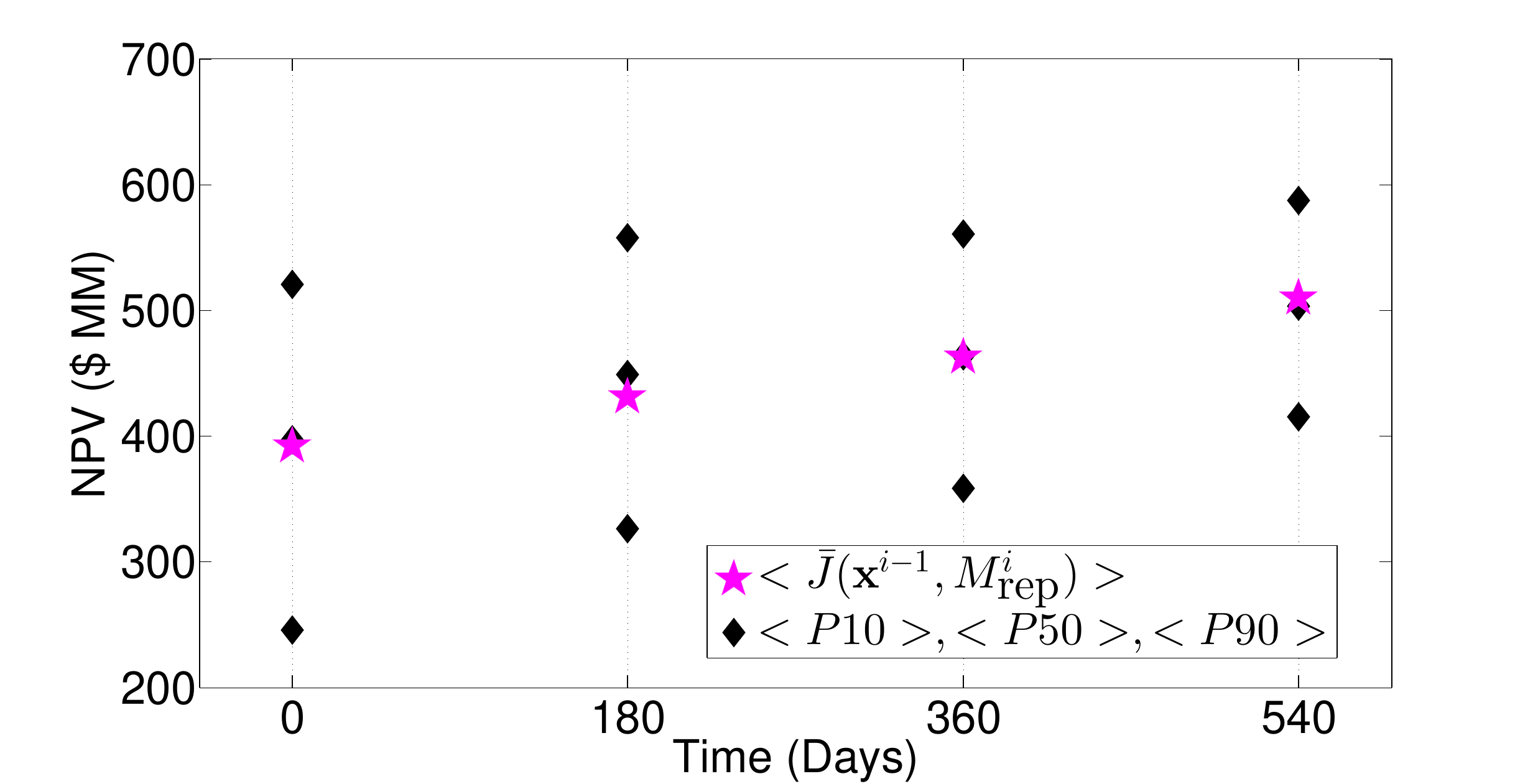}
\caption{Average values of P10, P50, P90 NPVs, over 25 runs along with the average true NPV, versus CLFD step.}
\label{fig:avg-P10P90}
\end{figure}

The average percentiles and the average expected NPV, for these 25 runs,
versus the CLFD step is shown in Fig~\ref{fig:avg-P10P90} and numerical values are provided in Table~\ref{tab:avg25}.
Each point in the figure (and table) is an average of 25 values (results in Fig.~\ref{fig:P10P90} shows values for run $j$=1), e.g.,
$<\bar{J}(\mathbf{x}^i,\mathbf{m}_{\text{true}})>$ = $\frac{1}{25}\sum_{j=1}^{25}J(\mathbf{x}^i_j, \mathbf{m}_{\text{true},j})$
where $\mathbf{x}^i_j$ is the optimal solution at CLFD step $i$ for run $j$, and $\mathbf{m}_{\text{true},j}$ is true-model $j$.
These results indicate a monotonic increase in all average values.
Further, the difference between the average P10 and the average P90 shrinks with CLFD step, although
this does not necessarily hold for each run.
These results are not specific to a particular true-model, and therefore they demonstrate the overall performance
of CLFD and present an assessment of project value under optimal development.

The author believes that this analysis
can help guide investment decisions to lease/buy an asset and in contract decisions between
independent oil companies (IOCs) and national oil companies (NOCs). In this example,
an investor may lease this reservoir for \$~246~million ($<P10>$ for OSV over prior), knowing that with CLFD, $<P10>$ at $t^{i=3}$ is \$~415~million.

\begin{table}
\centering \caption{
Numerical values of average P10, P50, P90 NPVs, over 25 runs along with the average true NPV, versus CLFD step.}\label{tab:avg25}
\begin{tabular}{lccccc}
  \hline
 & OSV-prior ($i$=0) & $t$=180 ($i$=1) & $t$=360 ($i$=2) & $t$=540 ($i$=3)\\
\noalign{\smallskip} \hline\noalign{\smallskip}
$<P90>$ &   521	& 558	& 561	& 588 \\
$<P50>$ &   397	& 449	& 464	& 504 \\
$<P10>$ & 	246	& 326	& 358	& 415 \\
$<\bar{J}(\mathbf{x}^i,\mathbf{m}_{\text{true}})>$ & 	392	& 431	& 463	& 510 \\
\hline\noalign{\smallskip}
\noalign{\smallskip}
\end{tabular}
\end{table}

The computational experiment of performing 25 CLFD applications took about 789~hours (33~days).
The total number of simulations is 9,333,900 which are performed through 30,136 calls to a computer cluster with access to maximum 400 cores.

\subsection{Discussion} \label{sec:disc}
The results for performance assessment showed a high range of possible outcomes
from CLFD runs.
The highly nonlinear nature of the optimization step
is perhaps the main cause of the wide range of variability in the ultimate outcome.
Although, the application of PSO component in the PSO-MADS algorithm,
provides some degree of global exploration, the optimization
may ultimately converge to a locally optimal solution (for the representative subset)
and there is no guarantee of finding the global optimum.
Results presented in \citet{onwunalu:10}
show that for a problem where the location
of a single well is optimized, the objective function surface is highly rough with many local optima.
The field development optimization problem in CLFD is significantly more complex, and a significant number of local optima may exist.
Another main cause to the wide range of final outcome,
is the sparsity in spatial data and the inherent uncertainty in model description after model calibration.
Despite all the complexity,
application of the full CLFD improved the true NPV in 96\% of cases, proving the effectiveness of the framework.

The statistical performance assessment procedure can be applied to assess the value of a subsurface reservoir
under optimal CLFD.
Decision making depends on the risk attitude and therefore, presenting a range
of possible outcomes is essential, e.g., for a risk-averse decision maker, the P10 estimate of the project value is more important
than the expected value.
By performing multiple CLFD runs, results similar to those presented in Table~\ref{tab:dJpr}
and Fig~\ref{fig:avg-P10P90} can be obtained.
These results present an assessment of
the project outcome when optimal reservoir management is to be applied.


The TruMAP statistical assessment procedure can be applied to test and compare the effect of different treatments.
Different geostatistical algorithms, optimization methods, representative model selection,
or model calibration methods can be applied in closed-loop optimization.
In order to compare two different treatments, multiple runs should be performed
and results can be analyzed in terms of a performance measure, e.g., average or median true-model improvement,
or difference in probability of increasing true-model NPV.

The statistical performance assessment procedure (TruMAP) and results represent an example of
massive computational experiments to quantify the effectiveness of a new algorithm. As \citet{monajemi:16} discussed,
with advancement of computational power and tools,
researchers are expected to present findings that involve such massive computations (order one million CPU-hours) as opposed to deducing conclusions
from a limited number of test cases.


\section{Conclusions and future work} \label{sec:CR}
In this work, a new implementation of closed-loop field development (CLFD) optimization is presented for models
described by multipoint geostatistics (MPS). CLFD involves three major steps: 1) optimization of
the full development plan for the entire project life, 2) drilling and completing $n_{\text{rig}}$ new wells and collecting new data,
and 3) updating models based on all data.
The optimization step is performed over multiple models through optimization with sample validation (OSV).
A two-step model calibration procedure is presented where spatial data (hard data at well locations)
are integrated through a MPS simulation method, and production data are integrated through a gradient-based RML procedure
using O-PCA parameterization.

The new TruMAP procedure is introduced for performance assessment of a closed-loop optimization algorithm.
Within this procedure, the new CLFD methodology is applied 25 times (for 25 different true-model cases),
and the CLFD performance was assessed in terms of true NPV improvement.
Results indicated that the true NPV in 64\%, 68\% and 96\% of cases increased
after one, two and three~steps, respectively.
These results showed that application of multi-step CLFD
is much more likely to increase true NPV than a single step CLFD,
particularly as richer data appears at later CLFD steps.
Further, the probability of increasing true-model NPV increases with CLFD step.
The statistical assessment results presented in this work, demonstrates the importance of
performing multiple runs
and considering several true-model cases to assess the performance
of any closed-loop optimization framework.


There are several directions to investigate in future.
Application of CLFD optimization to geothermal reservoir, groundwater remediation \citep{bayer:08,ghorbanidehno:17},
and CO$_2$-EOR problems, and unconventional resources \citep{shirangi:19} should also be investigated.
In this work, a fixed project life was specified in the optimization step of CLFD and the expected NPV was optimized.
It will be of interest to
assess CLFD performance when the rate of return is incorporated in optimization objective and
the economic project life is also optimized \citep{shirangi:17b}.
New implementations of closed-loop optimization should be statistically assessed by application of TruMAP procedure.


\section*{Acknowledgments}
This research was conducted at Stanford University between 2015 and 2017.
Access to the computational resources provided by the Stanford Center for Computational Earth \& Environmental Science (CEES)
is gratefully acknowledged.


\end{document}